\documentclass{amsproc}

\newtheorem{theorem}{Theorem}[section]
\newtheorem{thm}[theorem]{Theorem}

\newtheorem{lem}[theorem]{Lemma}
\newtheorem{cor}[theorem]{Corollary}
\newtheorem{prop}[theorem]{Proposition}

\theoremstyle{definition}
\newtheorem{definition}[theorem]{Definition}
\newtheorem{defn}[theorem]{Definition}
\newtheorem{example}[theorem]{Example}
\newtheorem{conj}{\bf Conjecture}

\theoremstyle{remark}

\newenvironment{pf}{\proof}{\endproof}

\numberwithin{equation}{section}

\begin{document}

\title[Symmetric Functions \& Representations of Quantum Affine Algebras]
{Symmetric Functions and Representations of Quantum Affine Algebras}

\author{Vyjayanthi Chari}
\address{Department of Mathematics\\ University of California
Riverside\\ Riverside, CA 92521}
\email{chari@newmath.ucr.edu}

\author{Michael Kleber}
\address{Department of Mathematics\\ Massachusetts Institute of
Technology\\ Cambridge, MA 02139}
\email{kleber@math.mit.edu}
\thanks{The second author was supported in part by an NSF
Mathematical Sciences Postdoctoral Research Fellowship.}

\date{21 November 2000}
\renewcommand{\subjclassname}{\textup{2000} Mathematics Subject Classification}

\begin{abstract}
We study connections between the ring of symmetric functions and the
characters of irreducible finite-dimensional representations of
quantum affine algebras.  We study two families of representations of
the symplectic and orthogonal Lie algebras.  One is defined via
combinatorial properties and is easy to calculate; the other is
closely related to the $q=1$ limit of the ``minimal affinization''
representations of quantum affine algebras.  We conjecture that the
two families are identical, and present supporting evidence and
examples.  In the special case of a highest weight that is a multiple
of a fundamental weight, this reduces to a conjecture of Kirillov and
Reshetikhin, recently proved by the first author.
\end{abstract}

\maketitle

\addtocounter{section}{-1}     
{\newcommand{\g}{\mathfrak{g}}

\section{Introduction}
\label{sec_intro}

In this paper we study connections between the ring of symmetric
functions and the characters of irreducible finite-dimensional
representations of quantum affine algebras.  We introduce the reader
to two families of representations of the classical finite-dimensional
simple Lie algebras $\g$, indexed by dominant integral weights
$\lambda$.  One family is defined via combinatorics and the ring of
symmetric functions, and is easy to describe and calculate. The other
family consists of ``minimal affinizations'' \cite{CP}, certain
representations of quantum affine algebras, regarded as
representations of the underlying finite-dimensional algebra. We
conjecture that these two families are identical and prove the
conjecture in certain cases.  In addition, we establish a number of
results which provide compelling evidence for the conjecture and also
illuminate the structure of the minimal affinizations of quantum
groups.

In Section~\ref{sec_sym} we define the representations
$W_{Sp}(\lambda)$ and $W_{O}(\lambda)$ of the symplectic and
orthogonal algebras, respectively.  They are described in terms of
their universal characters, which are elements of the ring of
symmetric functions.  They have the remarkable property that the map
taking the Schur function $s_\lambda$ to the character of
$W_{G}(\lambda)$ is an isomorphism of the ring of symmetric functions
($G=Sp$ or $O$).  This condition suffices to define the
representations $W_G(\lambda)$ completely.  In the special case when
$\lambda$ is a multiple of a fundamental weight (also called a
rectangle, from the shape of its Young diagram), the modules
$W_G(\lambda)$ were defined earlier by Kirillov and Reshetikhin
(\cite{KR}, see also~\cite{HKOTY}); the fact that this assignment
extends to a homomorphism of rings was proved in \cite{embed}.

In Section~\ref{sec_algebra} we define representations
$W_\g^{\text{aff}}(\lambda)$ of the loop algebra $L(\g)$.  They have
as quotients the $q=1$ specialization of the ``minimal affinization,''
a canonical representation of $U_q(\hat{\g})$ associated to each
$\lambda$.  In \cite{C} it was proved that when $\lambda$ is a
rectangle, they are in fact isomorphic and that
\begin{equation*}
W_\g^{\text{aff}}(\lambda)\cong W_G(\lambda),
\end{equation*}
a result that was conjectured in \cite{KR}.  In this paper, we
conjecture that this is true for all $\lambda$.  Most of the section
is devoted to proving results on the $\g$-module structure of the
modules $W_\g^{\text{aff}}(\lambda)$ which support this conjecture.
We isolate crucial properties that are known to be true for one family
and prove that the other family also satisfies them.
\begin{enumerate}
\item
The modules $W_G(\lambda)$ are defined for $so(n)$ and $sp(2n)$ for
all $n$.  Central to their definition is the fact that the direct sum
decomposition into irreducibles is independent of $n$, provided $n$ is
sufficiently large.  We show that the $W_\g^{\text{aff}}(\lambda)$ share
this property.
\item
The minimal affinization is distinguished from all other affinizations
of $\lambda$ by the property that if it contains a $\g$-highest weight
vector with weight $\mu$, then the root $\lambda-\mu$ cannot be
contained in a sub-root-lattice for a subalgebra isomorphic to some
$sl(r)$.  We show that $W_G(\lambda)$ and $W_\g^{\text{aff}}(\lambda)$
both have this property.
\item
Finally, in Section~\ref{sec_compute} we calculate examples.  We
restrict our attention to $\g=so(2n)$ for convenience; the proofs in
the other cases are similar.  For several types of weights $\lambda$,
we completely calculate $W_O(\lambda)$ and show explicitly that
$W_\g^{\text{aff}}(\lambda)$ is a submodule.
\end{enumerate}
We also mention other properties which we can prove for one of
$W_G(\lambda)$ or $W_\g^{\text{aff}}(\lambda)$ and which the other seems
empirically to share, though we cannot as yet  provide a proof.

The structure of finite-dimensional representations of quantum affine
algebras is very complicated.  Establishing the conjecture in this
paper would significantly expand our understanding of their algebraic
and combinatorial structure.  The work of \cite{KR} also conjectured a
formula for decomposing tensor products of representation associated
to rectangles, given by the fermionc formula. A generalization of that
formula beyond rectangles using the representations studied here would
be of considerable interest.

}
{\newcommand{\C}{\mathbf{C}}
\newcommand{\Z}{\mathbf{Z}}
\newcommand{\g}{\mathfrak{g}}
\renewcommand{\sl}{\mathfrak{sl}}
\renewcommand{\l}{\ell}
\newcommand{\pa}[1]{{\langle{#1}\rangle}}
\newcommand{\tensor}{\otimes}
\newcommand{\iso}{\cong}
\newcommand{\cont}{{\mathop{\mathrm{content}}\nolimits}}
%
%
\newcommand{\vdomsize}[2]{{\begin{picture}(#1,#2)\multiput(0,0)(0,#1){3}%
{\line(1,0){#1}}\multiput(0,0)(#1,0){2}{\line(0,1){#2}}\end{picture}}}
\newcommand{\hdomsize}[2]{{\begin{picture}(#2,#1)\multiput(0,0)(#1,0){3}%
{\line(0,1){#1}}\multiput(0,0)(0,#1){2}{\line(1,0){#2}}\end{picture}}}
\newcommand{\vdom}{{\vdomsize48}} \newcommand{\vdoms}{{\vdomsize36}}
\newcommand{\hdom}{{\hdomsize48}} \newcommand{\hdoms}{{\hdomsize36}}

%
%

\section{Symmetric Functions}
\label{sec_sym}

In this section we fix notation and recall the basic notions of the
ring $\Lambda$ of symmetric functions as a tool for handling
representations of the classical Lie algebras.  Our goal is the
definition of a certain subcategory of finite-dimensional
representations of the orthogonal and symplectic Lie algebras. This
subcategory is closed under taking direct sums and tensor products,
and it is generated as an abelian group by a family of modules
$W(\lambda)$, as $\lambda$ runs over all dominant integral highest
weights.  It has the remarkable property that the multiplicities in
the decomposition of tensor products are the Littlewood--Richardson
numbers.

\subsection{Some classical bases}
\label{ss_A}

We will work in the ring $\Lambda$ of formal symmetric functions in
countably many variables $(x_1,x_2,\ldots)$, and primarily follow the
notation of Macdonald~\cite{Macdonald}, to which we refer the reader
for proofs of fundamental facts.  Our emphasis is the dictionary which
translates between the combinatorics of symmetric functions and the
representation theory of the Lie algebra $\sl_n$.

The $k$th {\em complete} symmetric function $h_k$ is the sum of all
monomials of degree $k$ in the variables $(x_1,x_2,\ldots)$.
The $k$th {\em elementary} symmetric function $e_k$ is the sum of all
square-free monomials of degree $k$ in the variables $(x_1,x_2,\ldots)$.
These functions are clearly symmetric, {\em i.e.} invariant under all
permutations of the $x_i$.  Each of these sets is algebraically
independent, and $\Lambda$ is exactly the ring of polynomials in
either the $h$'s or the $e$'s.

We can specialize these function to polynomials by setting all
variables except for $x_1,\ldots,x_n$ to be zero, for any positive
integer $n$.  Now they are intimately familiar to representation
theorists: note that $h_1=e_1$ is the character of the fundamental
$n$-dimensional vector representation $V$ of the Lie group $GL(n)$ or
the Lie algebra $\sl_{n+1}$.  More generally, $h_k$ is the character
of $S^k(V)$, the $k$th symmetric power of the vector representation,
while $e_k$ is the character of $\bigwedge^k(V)$, the $k$th
alternating power; all of these are irreducible representations.

More precisely, recall that $\sl_{n+1}$ has $n$ {\em fundamental
weights\/} $\omega_1,\ldots,\omega_n$, and that its finite-dimensional
irreducible representations are indexed by {\em dominant} weights,
positive integer linear combinations of the fundamental weights.  We
will write $\lambda = a_1\omega_1 + \cdots + a_n\omega_n$ for such a
weight and $V(\lambda)$ for the associated representation; $\lambda$
is called the {\em highest weight}.  Then $h_k$ specializes to the
character of $V(k\omega_1)$, while $e_k$ specializes to the character
of $V(\omega_k)$.

The {\em Schur functions\/} fill out this picture.  For each dominant
weight $\lambda$, the Schur function $s_\lambda\in\Lambda$ is a
symmetric function in $(x_1,x_2,\ldots)$, and when we specialize to a
polynomial by setting $x_i=0$ for $i>n$, we get the character of the
representation $V(\lambda)$ of $\sl_{n+1}$.  Note that $n$ must be
large enough that $\lambda$ involves only the first $n$ fundamental
weights.  It is more traditional to index the Schur functions
$s_\lambda$ by {\em partitions}, integer sequences
$\lambda=\pa{\lambda_1,\ldots,\lambda_r}$ with
$\lambda_1\geq\cdots\geq\lambda_r\geq0$, or by their graphical
representation {\em Young diagrams}, in which $\lambda$ is depicted as
an array of $r$ rows of boxes, left-justified, with $\lambda_i$
boxes in the $i$th row.  Translation between a weight and a partition
is straightforward: the coefficient in $\lambda$ of $\omega_k$ is
$\lambda_k-\lambda_{k+1}$, or the number of columns of height exactly
$k$ in the Young diagram.  For example:
$$
\lambda=\omega_1+2\omega_2+\omega_3,\quad
\lambda=\pa{4,3,1},\quad
\mbox{Young diagram }=\,\,
{\setlength{\unitlength}{7pt}
\raisebox{-1.2\unitlength}{\begin{picture}(4,3)
\put(0,3){\line(1,0){4}}\put(0,2){\line(1,0){4}}
\put(0,1){\line(1,0){3}}\put(0,0){\line(1,0){1}}
\put(0,3){\line(0,-1){3}}\put(1,3){\line(0,-1){3}}\put(2,3){\line(0,-1){2}}
\put(3,3){\line(0,-1){2}}\put(4,3){\line(0,-1){1}}
\end{picture}}}
$$
Thus the finite-dimensional irreducible representations of $\sl_{n+1}$
are indexed by the Young diagrams with at most $n$ rows.

Since the ring of symmetric functions is the ring of polynomials in
the $h$'s, we must be able to write $s_\lambda$ as such a polynomial.
This is accomplished by the {\em Jacobi--Trudi} identity: 
$$
s_\lambda = \det\left(h_{\lambda_i-i+j}\right)_{i,j=1,\ldots,r},\quad
\mbox{ where }
\lambda=\pa{\lambda_1,\ldots,\lambda_r}.
$$
We can use this to write $s_\lambda$ in terms of the $e$'s as well
once we introduce the involutive ring automorphism
$\omega:\Lambda\to\Lambda$, defined by any one of the following:
$$
\omega(h_k)=e_k,\quad
\omega(e_k)=h_k,\quad
\omega(s_\lambda)=s_{\lambda'},
$$
where $\lambda'$ is the {\em conjugate} of $\lambda$, whose Young
diagram is obtained from that of $\lambda$ by reflecting it through
its main diagonal, exchanging rows and columns.  Note that
$h_k=s_\pa{k}$, whose Young diagram is a single row of length $k$, and
$e_k=s_\pa{1,\ldots,1}$, whose Young diagram is a single column of
height $k$.

The Schur functions $s_\lambda$ form a linear basis for $\Lambda$ as
$\lambda$ runs through all partitions (including the empty partition;
$s_\emptyset=1$ corresponds to the trivial representation, with
highest weight $0$).  We endow $\Lambda$ with an inner product
$\langle\,,\,\rangle$ by declaring this basis orthonormal.  The
multiplication in $\Lambda$ has as its structure constants the {\em
Littlewood--Richardson numbers:}
\begin{equation}
\label{defc}
c^\lambda_{\mu\nu}:=\langle s_\lambda,s_\mu s_\nu\rangle,\quad
\mbox{or equivalently,}\quad
s_\mu s_\nu = \sum_{\lambda} c^\lambda_{\mu\nu}\, s_\lambda.
\end{equation}
We write $|\lambda|$ for the number of boxes in the Young diagram of
$\lambda$, which is also the degree of the monomials in $s_\lambda$.
As a result,
\begin{equation}
\label{zero}
c^\lambda_{\mu\nu}=0 \quad\mbox{unless}\quad |\lambda|=|\mu|+|\nu|.
\end{equation}

Since Schur functions are characters, the Littlewood--Richardson
numbers describe the decomposition of a tensor product of
representations of $\sl_n$ into irreducibles: $c^\lambda_{\mu\nu}$ is
the multiplicity of $V(\lambda)$ in $V(\mu)\tensor V(\nu)$.  Thus the
Littlewood--Richardson numbers are nonnegative integers.

\subsection{Symplectic and orthogonal analogues}
\label{ss_BCD}

We now turn to the other classical Lie algebras, those of the
symplectic and orthogonal (of even or odd dimension) groups.  The work
of Koike and Terada~\cite{KT} showed that they too have ``universal
characters'' in $\Lambda$, which specialize to characters of their
representations in a way analogous to the Schur functions and $\sl_n$.

The irreducible finite-dimensional representations of the Lie groups
$SO(2n+1)$ and $Sp(2n)$ are once again indexed by highest weights
which are, as above, in bijection with Young diagrams with at most $n$
rows.  The irreducible representations of $O(2n)$ are indexed by the
same set, but when we restrict to $SO(2n)$, representations associated
to Young diagrams with exactly $n$ rows split into two irreducibles
(exchanged by the automorphism switching the two ``spin'' weights).
The starting point of~\cite{KT} is the observation that the characters
of these representations, just as in the case of $GL(n)$, are
``stable,'' in the sense that they are all specializations of their
$n\to\infty$ limit.

These ``stable limit'' or ``universal'' characters form two new bases
for $\Lambda$.  One consists of the characters $sp_\lambda$ coming
from the symplectic groups, and the other of the characters
$o_\lambda$ coming from the orthogonal groups (which give one stable
limit, independent of the parity of their rank).  We can specialize
$o_\lambda$ and $sp_\lambda$ to get characters of irreducible
orthogonal or symplectic representations of $SO(2n+1)$, $SO(2n)$ or
$Sp(2n)$ as long as $n$ is large enough; for our purposes we note that
$n$ is sufficiently large when the number of nonzero parts of
$\lambda$ is at most $n-1$ (or $n-2$ for $SO(2n)$).  When $n$ is too
small we get a character of a reducible representation; we refer
readers to the original paper or a well-written summary, like
Appendix~A of~\cite{FH}, for this level of details of the
specialization homomorphisms.

These two new bases have their own structure constants, the symplectic
and orthogonal analogues of the Littlewood--Richardson numbers.  The
following remarkable fact deserves wider recognition.

\begin{theorem}
\label{thm_spo}
There is a collection of nonnegative integers $d^\lambda_{\mu\nu}$
such that
$$
sp_\mu\,sp_\nu = \sum_{\lambda} d^\lambda_{\mu\nu}\, sp_\lambda
\quad\mbox{ and }\quad
o_\mu\,o_\nu = \sum_{\lambda} d^\lambda_{\mu\nu}\, o_\lambda.
$$
That is, the structure constants of the symplectic and orthogonal
bases are the same.
\end{theorem}

The equality is shown in~\cite{KT} in terms of symmetric functions; it
also follows easily from crystal base theory.  The
$d^\lambda_{\mu\nu}$ are certainly nonnegative integers, as they count
the multiplicity of $V_G(\lambda)$ in $V_G(\mu)\tensor V_G(\nu)$.  Here the
$V_G$ can denote representations of any one of $G=O(2n+1)$, $Sp(2n)$, or
$O(2n)$, so long as we require $n$ to be sufficiently large --- in
particular, larger than the sum of the numbers of rows in $\mu$ and
$\nu$.  For example, 
$$
\begin{array}{ll}
V(\omega_1)^{\tensor2} \iso V(2\omega_1)\oplus V(\omega_2) &
  \mbox{in $GL(n)$ or $SL(n+1)$}, n\geq 2, \\[4pt]
V_G(\omega_1)^{\tensor2} \iso V_G(2\omega_1)\oplus V_G(\omega_2) \oplus \C&
  \mbox{in $SO(2n+1)$, $Sp(2n)$, or $SO(2n)$},  n>2.
\end{array}
$$
Here $\C$ denotes the trivial representation $V_G(0)$.  

In some sense the $d^\lambda_{\mu\nu}$ are a deformation of the
Littlewood--Richardson numbers.  Analogous to \eqref{zero}, we have
\begin{equation}\label{dmunu}
d^\lambda_{\mu\nu}=0 \quad\mbox{unless}\quad
|\lambda|=|\mu|+|\nu|-2k,\,\, k\in\Z_{\geq0}.
\end{equation}
Recall that for $c^\lambda_{\mu\nu}$ we demand this with $k=0$.  Moreover,
\begin{equation}\label{cd}
c^\lambda_{\mu\nu} = d^\lambda_{\mu\nu}
\quad\mbox{when}\quad |\lambda|=|\mu|+|\nu|.
\end{equation}
In other words, moving from the general linear
to the symplectic or orthogonal groups only adds new pieces to the
decomposition of tensor products, and all the new pieces are
lower-order terms.  

From the representation theory point of view, $\lambda$ is a weight
and can be written as a linear combination of the fundamental roots
$\alpha_1,\ldots,\alpha_n$.  As long as $k<n$ (or $n-1$, for type
$D_n$), the coefficient of $\alpha_k$ is the number of boxes in the
top $k$ rows of the Young diagram of $\lambda$, so $|\lambda|$ is the
coefficient of $\alpha_k$ for $k$ greater than the number of rows of
$\lambda$.  (The coefficient of $\alpha_n$ may differ from $|\lambda|$ 
but only by a factor of two.)  Therefore the ``extra pieces'' of
$V(\mu)\tensor V(\nu)$ corresponding to nonzero $d^\lambda_{\mu\nu}$
with $|\lambda|<|\mu|+|\nu|$ can be identified by the fact that the
weight $\mu+\nu-\lambda$ is supported on the ``spin'' or ``long''
root $\alpha_n$ (and for $SO(2n)$ on $\alpha_{n-1}$ also).

\subsection{A new family of representations}
\label{ss_defW}

We are now ready to define the symplectic and orthogonal families of
representations $W(\lambda)$.  The stable limit characters of the
$W$'s will form another pair of new bases of the ring $\Lambda$, and
the representations are completely characterized by the property that
the structure constants of these new bases are the classical
Littlewood--Richardson numbers $c^\lambda_{\mu\nu}$.

Consider the natural inclusions $SO(2n+1)\subseteq GL(2n+1)$,
$Sp(2n)\subseteq GL(2n)$, and $SO(2n)\subseteq GL(2n)$.  In each case
the inclusion gives rise to a restriction map which takes any
representation of the general linear group and views it as a module
over the symplectic or orthogonal subgroup.  If we pick one of these
$G\subseteq GL(n)$ and an irreducible $GL(n)$ module $V(\lambda)$, its
restriction $V(\lambda)|_G$ will in general be reducible, and the
decomposition into symplectic or orthogonal irreducibles is
independent of $n$ as long as $\lambda$ has at most $n$ rows.

Taking characters translates this decomposition into the question of
writing the Schur functions in the $sp_\lambda$ or $o_\lambda$ bases;
the coefficients will be the multiplicities, so will certainly be
nonnegative integers.  The ``branching rules'' were known to
Littlewood; they were written in the context of symmetric functions by
Koike and Terada.

\begin{theorem}
\label{thm_cit}
The following summations are over all partitions $\mu$.
\begin{enumerate}
\item
$s_\lambda = \sum_\mu \left(\sum_{\nu\in Y_\vdoms}
c^\lambda_{\mu\nu}\right) sp_\mu$, where $\nu\in Y_\vdom$ if the Young
diagram of $\nu$ has only even-height columns, {\em i.e.}~can be tiled
by vertical dominos.
\item
$s_\lambda = \sum_\mu \left(\sum_{\nu\in Y_\hdoms}
c^\lambda_{\mu\nu}\right) o_\mu$, where $\nu\in Y_\hdom$ if the Young
diagram of $\nu$ has only even-length rows, {\em i.e.}~can be tiled
by horizontal dominos.
\end{enumerate}
\end{theorem}

The names $Y_\vdom$ and $Y_\hdom$ are mnemonic, but they may be more
familiar by other names.  Viewed as sets of partitions, $\nu\in Y_\vdom$
means all parts of $\nu$ occur with even multiplicity, while $\nu\in
Y_\hdom$ means $\nu$ has exclusively even parts.  Viewed in terms of
weights, $\nu\in Y_\vdom$ means $\nu$ is in the span of the
even fundamental weights $\omega_{2i}$, while $\nu\in Y_\hdom$ means
the coefficient of $\omega_i$ in $\nu$ is even, for all $i$.

Now recall the remarkable fact from Theorem~\ref{thm_spo} that the
$sp_\lambda$ and $o_\lambda$ bases have the same structure constants.
The linear maps $sp_\lambda\mapsto o_\lambda$ and $o_\lambda\mapsto
sp_\lambda$ are therefore ring isomorphisms.  The representations we
are related to the irreducible representations $V(\lambda)$ by these
isomorphisms.

\begin{definition}
\label{def_W}
We define two families of reducible representations, by giving their
direct sum decomposition (with multiplicities) into irreducibles
$V_{Sp}$ or $V_{O}$ of the symplectic or orthogonal groups, respectively:
$$
W_{Sp}(\lambda) := \sum_\mu \left(\sum_{\nu\in Y_\hdoms}
c^\lambda_{\mu\nu}\right) V_{Sp}(\mu)
\qquad
W_{O}(\lambda) := \sum_\mu \left(\sum_{\nu\in Y_\vdoms}
c^\lambda_{\mu\nu}\right) V_{O}(\mu)
$$
\end{definition}

One needs to compare Theorem~\ref{thm_cit} with Definition~\ref{def_W}
carefully to distinguish the irreducible $GL$-module $V(\lambda)$ from
these new $W_G(\lambda)$: the two differ typographically only by
exchanging $Y_\vdom$ with $Y_\hdom$.  To avoid confusion, we place the
four possibilities in one table to highlight their relations:
\begin{equation}
\label{VW}
\renewcommand{\arraystretch}{1.5}
\begin{array}{l|cc}
\sum c^\lambda_{\mu\nu} &  V_{Sp}(\mu)      &  V_{O}(\mu) \\
\hline
     \nu\in Y_\vdom     &  V(\lambda)       &  W_{O}(\lambda)\\
     \nu\in Y_\hdom     &  W_{Sp}(\lambda)  &  V(\lambda)
\end{array}
\end{equation}

Note that the decomposition of $W_G(\lambda)$ is independent of $n$ as
long as $\lambda$ has fewer than $n$ rows.  Therefore we can speak of
the $W_G(\lambda)$ giving rise to stable limit characters, just as the
$V(\lambda)$ give rise to the Schur functions as described in
Section~\ref{ss_A}.  In light of the similarity in their respective
definitions, it is immediate that the stable limit characters of
$W_{Sp}(\lambda)$ and $W_{O}(\lambda)$ are the images of $s_\lambda$
under the ring isomorphisms $o_\mu\mapsto sp_\mu$ and $sp_\mu\mapsto
o_\mu$, respectively.  Because their characters are images of the
Schur functions, we have a family of reducible representations of the
symplectic or orthogonal groups whose tensor products decompose into a
direct sum of family members according to the classical
Littlewood--Richardson numbers:
\begin{equation}
\label{Wlr}
W(\mu) \tensor W(\nu) \iso \sum_\lambda c^\lambda_{\mu\nu}\,W(\lambda), \quad
\mbox{where }\, W = W_{Sp} \mbox{ or } W_{O}.
\end{equation}
In~\cite{embed} we showed that the $W(\lambda)$ are completely
characterized by this property:

\begin{theorem}
Let $\{X(\lambda)\}$ be a family of representations of the symplectic
or orthogonal groups $G$, indexed by all partitions $\lambda$ and
given in terms of their irreducible decompositions $X(\lambda) \iso
V_G(\lambda) \oplus
\sum_{\mu<\lambda} m_{\lambda\mu} V_G(\mu)$, for some nonnegative
integers $m_{\lambda\mu}$.  Suppose their tensor products decompose
into direct sums as $X(\mu) \tensor X(\nu) \iso \sum_\lambda
c^\lambda_{\mu\nu}\, X(\lambda)$.  Then there are only two possibilities:
\begin{enumerate}
\renewcommand{\theenumi}{\alph{enumi}}\renewcommand{\labelenumi}{(\theenumi)}
\item
$\{X(\lambda)\}$ is the family of $GL$-irreducible representations
$\{V(\lambda)\}$; or
\item
$\{X(\lambda)\}$ is the family $\{W_G(\lambda)\}$ defined above.
\end{enumerate}
\end{theorem}

We have defined $W_G(\lambda)$ to be a representation of $G$, but for
most of the remainder of this paper we will be primarily concerned
with Lie algebras (and their loop algebras and quantum deformations).
We therefore allow $W_G(\lambda)$ to denote a representation of the
finite-dimensional Lie algebra $\g$ as well.  Since we defined the
representation in terms of its direct sum decomposition into
irreducibles, nothing new is introduced by this convenience.

\subsection{Combinatorics of $\mathbf{W_G}$}
\label{ss_combW}

We will need the ability to compute $W_G(\lambda)$ explicitly.
The definition in terms of Littlewood--Richardson numbers can be
restated in terms of the {\em skewing} operation
$s_\nu^\perp:\Lambda\to\Lambda$, the adjoint to multiplication by
$s_\nu$.  The skew Schur function $s_{\lambda/\nu}$ is defined as
$s_\nu^\perp s_\lambda$.  The skew Young diagram for $\lambda/\nu$ is
represented by the Young diagram for $\lambda$ with the boxes for the
Young diagram of $\nu$ removed from the upper-left corner; this
notation relies on the fact that $s_{\lambda/\nu}=0$ unless
$\nu\subseteq\lambda$, where $\subseteq$ denotes containment of Young
diagrams.  We mention that $s_{\lambda/\nu}$ also has a Jacobi--Trudi
expansion:
$$
s_{\lambda/\nu} = \det(h_{\lambda_i-\nu_j-i+j})_{i,j=1,\ldots,r}.
$$

The character of $W_{O}(\lambda)$ can now be described as the image of
$\sum_{\nu\in Y_\vdoms} s_{\lambda/\nu}$ under the linear map $s_\mu
\mapsto o_\mu$, and similarly for $W_{Sp}(\lambda)$ with $\nu\in
Y_\hdom$ and $s_\mu \mapsto sp_\mu$.  Since we will focus on the
decomposition into irreducibles, we just need to expand $\sum_{\nu\in
Y_\vdoms} s_{\lambda/\nu}$ (or $\nu\in Y_\hdom$) in the Schur basis.
To calculate this character for a given $\lambda$, we will give a
combinatorial algorithm to expand each $s_{\lambda/\nu}$ as a sum of
Schur functions.  Our presentation is heavily abridged, to say the
least; we refer the reader to Chapter~5 of~\cite{Fulton} for a
rational and justified development.

\begin{defn}
\begin{enumerate}
\renewcommand{\theenumi}{\alph{enumi}}\renewcommand{\labelenumi}{(\theenumi)}
\item
A {\em semi-standard Young tableau} of shape $\lambda/\nu$ is a
filling of the boxes of the skew Young diagram of $\lambda/\nu$ with
nonnegative integers such that the entries {\em strictly increase}
reading down any column and {\em weakly increase} reading across any
row.
\item
The {\em reverse row word} of a tableau $T$ is the sequence of
integer entries of $T$ as you read each row from right to left,
beginning with the top row and ending with the bottom.
\item
A sequence of integers is a {\em ballot sequence} if the number of
occurrences of $i+1$ in the first $k$ terms is never greater than the
number of occurrences of $i$, for every integer $i$ and $k$.
\item
The {\em content} of a tableau $T$ is $\pa{n_1(T),n_2(T),\ldots}$,
where $n_i(T)$ is the number of occurrences of $i$ in $T$.
\end{enumerate}
\end{defn}

\begin{prop} $$
s_{\lambda/\nu} = \sum_T s_{\cont(T)},
$$
where $T$ ranges over all semi-standard Young tableaux of shape
$\lambda/\nu$ such that the reverse row word of $T$ is a ballot
sequence.  Since the reverse row word of $T$ is a ballot sequence,
$\cont(T)$ is a partition.
\end{prop}

To illustrate, let us compute $W_O(\lambda)$ for
$\lambda=\omega_1+\omega_2+\omega_3$. The only partitions $\nu\in
Y_\vdom$ contained in $\lambda$ are the empty partition, $\pa{1,1}$,
and $\pa{2,2}$.  The six semi-standard  Young tableaux $T$ whose
reverse row words are ballot sequences are:
$$
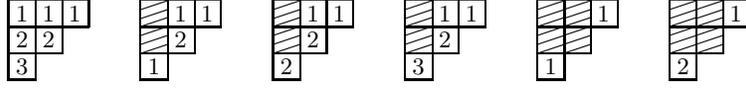
\setlength{\unitlength}{10pt}
\newcommand{\Yd}[1]{\begin{picture}(3,3)\put(0,0){\framebox(1,3){}}
\put(0,1){\framebox(2,2){}}\put(0,2){\framebox(3,1){}}#1\end{picture}}
\newcommand{\scribble}[2]{\multiput(#1,#2)(0,.33){6}{\line(3,1){1}}}
\newcommand{\tabl}[3]{\put(#1,#2){\makebox(1,1){\small #3}}}
\Yd{\tabl221\tabl121\tabl021\tabl112\tabl012\tabl003}
\qquad
\Yd{\scribble01\tabl221\tabl121\tabl112\tabl001}
\qquad
\Yd{\scribble01\tabl221\tabl121\tabl112\tabl002}
\qquad
\Yd{\scribble01\tabl221\tabl121\tabl112\tabl003}
\qquad
\Yd{\scribble01\scribble11\tabl221\tabl001}
\qquad
\Yd{\scribble01\scribble11\tabl221\tabl002}
$$
To convert $\cont(T)$ into a weight, recall that the coefficient of
$\omega_i$ is the number of occurrences of $i$ minus the number of
occurrences of $i+1$.  So $W_O(\omega_1+ \omega_2+ \omega_3)$
decomposes as a sum of six components, respectively:
$$
 V_O(\omega_1+ \omega_2+ \omega_3) \oplus V_O(2\omega_1+ \omega_2) \oplus
 V_O(2\omega_2) \oplus V_O(\omega_1+ \omega_3) \oplus V_O(2\omega_1) 
 \oplus V_O(\omega_2).
$$

\begin{prop}\label{motive}
Immediate consequences of the combinatorics of $W_G(\lambda)$:
\begin{enumerate}
\item
If $V_G(\mu)$ appears in $W_G(\lambda)$ then $\mu\subseteq\lambda$, where
$\subseteq$ denotes containment of Young diagrams.  In particular, if
$\mu\neq\lambda$ then $|\mu|<|\lambda|$.  As we observed at the end of
Section~\ref{ss_BCD}, this is equivalent to saying that the root
$\lambda-\mu$ is supported on the root $\alpha_n$ (and for $SO(2n)$
on $\alpha_{n-1}$ also).
\item
The trivial representation $V_G(0)=\C$ appears in $W_O(\lambda)$ if and
only if $\lambda\in Y_\vdom$, and in $W_{Sp}(\lambda)$ if and only if
$\lambda\in Y_\hdom$.  Moreover, if it does appear, it has
multiplicity $1$.  This follows because $s_{\lambda/\nu}$ is
homogeneous of degree $|\lambda|-|\nu|$, which by the previous point
is zero if and only if $\nu=\lambda$ {\rm(}when $s_{\lambda/\nu}=1${\rm)}.
\end{enumerate}
\end{prop}

The first consequence will be a significant motivation for the
conjecture in the next section, where another family of
representations is uniquely characterized by the similar property that
$\lambda-\mu$ cannot be contained in a sublattice corresponding to a
type $A$ subalgebra.  We do not know of a way to see the second
consequence from the representation theory side.
}
{
%

\newcommand{\thmref}[1]{Theorem~\ref{#1}}
\newcommand{\secref}[1]{Sect.~\ref{#1}}
\newcommand{\lemref}[1]{Lemma~\ref{#1}}
\newcommand{\propref}[1]{Proposition~\ref{#1}}
\newcommand{\corref}[1]{Corollary~\ref{#1}}
\newcommand{\remref}[1]{Remark~\ref{#1}}
\newcommand{\nc}{\newcommand}
\newcommand{\rnc}{\renewcommand}
\nc{\cal}{\mathcal}
\nc{\goth}{\mathfrak}
\rnc{\bold}{\mathbf}
\renewcommand{\frak}{\mathfrak}
\renewcommand{\Bbb}{\mathbb}

\nc{\Cal}{\mathcal}
\nc{\Xp}[1]{X^+(#1)}
\nc{\Xm}[1]{X^-(#1)}
\nc{\on}{\operatorname}
\nc{\ch}{\mbox{ch}}
\nc{\Z}{{\bold Z}}
\nc{\J}{{\cal J}}
\nc{\C}{{\bold C}}
\nc{\Q}{{\bold Q}}
\renewcommand{\P}{{\cal P}}
\nc{\N}{{\Bbb N}}
\nc\beq{\begin{equation}}
\nc\enq{\end{equation}}
\nc\lan{\langle}
\nc\ran{\rangle}
\nc\bsl{\backslash}
\nc\mto{\mapsto}
\nc\lra{\leftrightarrow}
\nc\hra{\hookrightarrow}
\nc\sm{\smallmatrix}
\nc\esm{\endsmallmatrix}
\nc\sub{\subset}
\nc\ti{\tilde}
\nc\nl{\newline}
\nc\fra{\frac}
\nc\und{\underline}
\nc\ov{\overline}
\nc\ot{\otimes}
\nc\bbq{\bar{\bq}_l}
\nc\bcc{\thickfracwithdelims[]\thickness0}
\nc\ad{\text{\rm ad}}
\nc\Ad{\text{\rm Ad}}
\nc\Hom{\text{\rm Hom}}
\nc\End{\text{\rm End}}
\nc\Ind{\text{\rm Ind}}
\nc\Res{\text{\rm Res}}
\nc\Ker{\text{\rm Ker}}
\rnc\Im{\text{Im}}
\nc\sgn{\text{\rm sgn}}
\nc\tr{\text{\rm tr}}
\nc\Tr{\text{\rm Tr}}
\nc\supp{\text{\rm supp}}
\nc\card{\text{\rm card}}
\nc\bst{{}^\bigstar\!}
\nc\he{\heartsuit}
\nc\clu{\clubsuit}
\nc\spa{\spadesuit}
\nc\di{\diamond}

\nc\al{\alpha}
\nc\bet{\beta}
\nc\ga{\gamma}
\nc\de{\delta}
\nc\ep{\epsilon}
\nc\io{\iota}
\nc\om{\omega}
\nc\si{\sigma}
\rnc\th{\theta}
\nc\ka{\kappa}
\nc\la{\lambda}
\nc\ze{\zeta}

\nc\vp{\varpi}
\nc\vt{\vartheta}
\nc\vr{\varrho}

\nc\Ga{\Gamma}
\nc\De{\Delta}
\nc\Om{\Omega}
\nc\Si{\Sigma}
\nc\Th{\Theta}
\nc\La{\Lambda}
\nc\boa{\bold a}
\nc\bob{\bold b}
\nc\boc{\bold c}
\nc\bod{\bold d}
\nc\boe{\bold e}
\nc\bof{\bold f}
\nc\bog{\bold g}
\nc\boh{\bold h}
\nc\boi{\bold i}
\nc\boj{\bold j}
\nc\bok{\bold k}
\nc\bol{\bold l}
\nc\bom{\bold m}
\nc\bon{\bold n}
\nc\boo{\bold o}
\nc\bop{\bold p}
\nc\boq{\bold q}
\nc\bor{\bold r}
\nc\bos{\bold s}
\nc\bou{\bold u}
\nc\bov{\bold v}
\nc\bow{\bold w}
\nc\boz{\bold z}

\nc\ba{\bold A}
\nc\bb{\bold B}
\nc\bc{\bold C}
\nc\bd{\bold D}
\nc\be{\bold E}
\nc\bg{\bold G}
\nc\bh{\bold h}
\nc\bH{\bold H}

\nc\bi{\bold I}
\nc\bj{\bold J}
\nc\bk{\bold K}
\nc\bl{\bold L}
\nc\bm{\bold M}
\nc\bn{\bold N}
\nc\bo{\bold O}
\nc\bp{\bold P}
\nc\bq{\bold Q}
\nc\br{\bold R}
\nc\bs{\bold S}
\nc\bt{\bold T}
\nc\bu{\bold U}
\nc\bv{\bold v}
\nc\bV{\bold V}

\nc\bw{\bold W}
\nc\bz{\bold Z}
\nc\bx{\bold x}
\nc\bX{\bold X}
\nc\blambda{{\mbox{\boldmath $\Lambda$}}}
\nc\bpi{{\mbox{\boldmath $\pi$}}}

\nc\e[1]{E_{#1}}
\nc\ei[1]{E_{\delta - \alpha_{#1}}}
\nc\esi[1]{E_{s \delta - \alpha_{#1}}}
\nc\eri[1]{E_{r \delta - \alpha_{#1}}}
\nc\ed[2][]{E_{#1 \delta,#2}}
\nc\ekd[1]{E_{k \delta,#1}}
\nc\emd[1]{E_{m \delta,#1}}
\nc\erd[1]{E_{r \delta,#1}}

\nc\ef[1]{F_{#1}}
\nc\efi[1]{F_{\delta - \alpha_{#1}}}
\nc\efsi[1]{F_{s \delta - \alpha_{#1}}}
\nc\efri[1]{F_{r \delta - \alpha_{#1}}}
\nc\efd[2][]{F_{#1 \delta,#2}}
\nc\efkd[1]{F_{k \delta,#1}}
\nc\efmd[1]{F_{m \delta,#1}}
\nc\efrd[1]{F_{r \delta,#1}}
\nc{\ug}{\bu^{fin}}

\nc\fa{\frak a}
\nc\fb{\frak b}
\nc\fc{\frak c}
\nc\fd{\frak d}
\nc\fe{\frak e}
\nc\ff{\frak f}
\nc\fg{\frak g}
\nc\fh{\frak h}
\nc\fj{\frak j}
\nc\fk{\frak k}
\nc\fl{\frak l}
\nc\fm{\frak m}
\nc\fn{\frak n}
\nc\fo{\frak o}
\nc\fp{\frak p}
\nc\fq{\frak q}
\nc\fr{\frak r}
\nc\fs{\frak s}
\nc\ft{\frak t}
\nc\fu{\frak u}
\nc\fv{\frak v}
\nc\fz{\frak z}
\nc\fx{\frak x}
\nc\fy{\frak y}

\nc\fA{\frak A}
\nc\fB{\frak B}
\nc\fC{\frak C}
\nc\fD{\frak D}
\nc\fE{\frak E}
\nc\fF{\frak F}
\nc\fG{\frak G}
\nc\fH{\frak H}
\nc\fJ{\frak J}
\nc\fK{\frak K}
\nc\fL{\frak L}
\nc\fM{\frak M}
\nc\fN{\frak N}
\nc\fO{\frak O}
\nc\fP{\frak P}
\nc\fQ{\frak Q}
\nc\fR{\frak R}
\nc\fS{\frak S}
\nc\fT{\frak T}
\nc\fU{\frak U}
\nc\fV{\frak V}
\nc\fZ{\frak Z}
\nc\fX{\frak X}
\nc\fY{\frak Y}
\nc\tfi{\ti{\Phi}}
\nc\bF{\bold F}

\nc\ua{\bold U_\A}

\nc\qinti[1]{[#1]_i}
\nc\q[1]{[#1]_q}
\nc\xpm[2]{E_{#2 \delta \pm \alpha_#1}}  
\nc\xmp[2]{E_{#2 \delta \mp \alpha_#1}}
\nc\xp[2]{E_{#2 \delta + \alpha_{#1}}}
\nc\xm[2]{E_{#2 \delta - \alpha_{#1}}}
\nc\hik{\ed{k}{i}}
\nc\hjl{\ed{l}{j}}
\nc\qcoeff[3]{\left[ \begin{smallmatrix} {#1}& \\ {#2}& \end{smallmatrix}
\negthickspace \right]_{#3}}
\nc\qi{q}
\nc\qj{q}

\nc\ufdm{{_\ca\bu}_{\rm fd}^{\le 0}}


\nc\isom{\cong} 

\nc{\pone}{{\Bbb C}{\Bbb P}^1}
\nc{\pa}{\partial}
\def\H{\cal H}
\def\L{\cal L}
\nc{\F}{{\cal F}}
\nc{\Sym}{{\goth S}}
\nc{\A}{{\cal A}}
\nc{\arr}{\rightarrow}
\nc{\larr}{\longrightarrow}

\nc{\ri}{\rangle}
\nc{\lef}{\langle}
\nc{\W}{{\cal W}}
\nc{\uqatwoatone}{{U_{q,1}}(\su)}
\nc{\uqtwo}{U_q(\goth{sl}_2)}
\nc{\dij}{\delta_{ij}}
\nc{\divei}{E_{\alpha_i}^{(n)}}
\nc{\divfi}{F_{\alpha_i}^{(n)}}
\nc{\Lzero}{\Lambda_0}
\nc{\Lone}{\Lambda_1}
\nc{\ve}{\varepsilon}
\nc{\phioneminusi}{\Phi^{(1-i,i)}}
\nc{\phioneminusistar}{\Phi^{* (1-i,i)}}
\nc{\phii}{\Phi^{(i,1-i)}}
\nc{\Li}{\Lambda_i}
\nc{\Loneminusi}{\Lambda_{1-i}}
\nc{\vtimesz}{v_\ve \otimes z^m}

\nc{\asltwo}{\widehat{\goth{sl}_2}}
\nc\eh{\frak h^e}  
\nc\loopg{L(\frak g)}  
\nc\eloopg{L^e(\frak g)} 
\nc\ebu{\bu^e} 
\nc\loopa{L(\frak a)}  

\nc\teb{\tilde E_\boc}
\nc\tebp{\tilde E_{\boc'}}

\section {Loop algebras and their representations} 
\label{sec_algebra}

From now on, $\frak g$ denotes a symplectic or orthogonal algebra of rank $n$ and $\loopg =\frak g\otimes\bc[t,t^{-1}]$ the corresponding loop algebra with the obvious Lie algebra structure.  
In this section, we construct a family of finite-dimensional indecomposable representations $W^{\text{aff}}(\lambda)$ of $\loopg$ and conjecture that these modules are isomorphic as  $\frak g$-modules to the modules $W(\lambda)$ of the previous section.

\subsection{Loop algebras}
 Fix a Cartan subalgebra $\frak h$ of $\frak g$ and a set of simple roots, $\{\alpha_1,\alpha_2,\cdots \alpha_n\}$, we shall assume that the simple roots are numbered as in \cite{H}.  Thus in the case of $sp(2n)$ (resp. $so(2n+1)$), we assume that $\alpha_n$ is the long (resp. short) root, while  in the case of $so(2n)$, we assume that $n-1$ and $n$ are the spin nodes. In all cases, we assume in addition that 
\begin{equation*} a_{i,i\pm 1}=-1, \ 1\le i\le n-2,\ \ a_{ij} =0, \ j\ne i\pm 1.\end{equation*} 
Set $J=\{\alpha_1,\cdots, \alpha_{n-2}, \alpha_{n-1}\}$.  In the case of $so(2n)$,
we shall also need the subset $\overline{J}= \{\alpha_1,\cdots, \alpha_{n-2}, \alpha_{n}\}$. 

Let $R^+$ be the set of positive roots and let  $Q^+$ (resp. $P^+$) be the integral root  (resp. weight) lattice respectively. Let $\omega_1,\cdots ,\omega_n$, be a set of fundamental weights in $P^+$. Given any subset $ J'$ of $\alpha_1,\alpha_2,\cdots \alpha_n$, let $R^+(J')$ be the subset of $R^+$ spanned by elements of $J'$. We define $Q^+(J')$ etc. in the obvious way. 

 For each $\alpha\in R^+$, fix nonzero elements $x_\alpha^\pm\in\frak g_{\pm\alpha}$ and  $h_\alpha\in\frak h$ satisfying,
\begin{equation*} [h_\alpha, x_{\alpha}^\pm] =\pm\alpha(h)x_{\alpha}^\pm, \ \ [x_\alpha^+,x_\alpha^-]=h_\alpha.\end{equation*}
 If $\alpha,\beta\in R^+$ is such that $\alpha+\beta\in R^+$ or $\alpha-\beta\in R^+$ we shall assume that \begin{equation*}[x_\alpha^\pm, x_\beta^\pm]=x_{\alpha+\beta}^\pm,\ \ \ \ [x_\alpha^\pm, x_\beta^\mp]=x_{\alpha-\beta}^\pm.\end{equation*}
Strictly speaking, the preceding statement is true only up to some nonzero scalar multiples, but for our purposes there is no loss in assuming, to simplify notation, that these scalars are all one.

 We shall need the following subalgebras of $\frak g$:
 \begin{align*}\frak k =\frak h\oplus_{\alpha\in R^+(J)}\frak g_\alpha,\ \  &\ \ \frak p=\frak k\oplus\frak u^+,\\  \frak u^\pm=\oplus_{\alpha\in R^+\backslash R^+(J)}\frak g_{\pm\alpha},\ \ &\ \ \frak n^\pm=\oplus_{\alpha\in R^+} \frak g_{\pm\alpha}.\end{align*} 
Then $\frak u^+$ is an ideal in $\frak p$ and we have a  homomorphism of Lie algebras $L(\frak p)\to L(\frak k)$ with kernel $L(\frak u^+)$. The algebra $\frak k$ is reductive, and hence 
 given any $\lambda\in P^+$ there exists an  irreducible $\frak k$-module $V_{\frak k}(\lambda)$ with highest weight $\lambda$ and highest weight vector $v_\lambda$. This can clearly be regarded as a module for $L(\frak k)$ by composing with the evaluation homomorphism $ev:L(\frak k)\to \frak k$ that sends $x\otimes t^n\to x$ and hence also as  a module for $L(\frak p)$.

\subsection{Construction of ${\mathbf W}^{\text{aff}}$} Consider the induced module
\begin{equation*}
 {\text{Ind}}_{L(\frak p)}^{L(\frak g)}V_{\frak k}(\lambda)= \bu(\loopg)\otimes_{L(\frak p)} V_{\frak k}(\lambda),\end{equation*}
where for any Lie algebra $\frak a$, we let $\bu(\frak a)$ be the universal enveloping algebra of $\frak a$. 
 Since
\begin{equation*}(x_{\alpha}^+\otimes t^n)(1\otimes v_\lambda) =0, \ \ \forall\  \alpha\in R^+, n\in\bz,\end{equation*}
an elementary application of the Poincare--Birkhoff--Witt  theorem shows that,
\begin{equation*} {\text{Ind}}_{L(\frak p)}^{L(\frak g)}V_{\frak k}(\lambda) =\bu(L(\frak n^-))(1\otimes v_\lambda)=\frak\bu(L(\frak u^-))\otimes V_{\frak k}(\lambda).\end{equation*} 
 In particular, the subspace  \begin{equation*}\left({\text{Ind}}_{L(\frak p)}^{L(\frak g)}V_{\frak k}(\lambda)\right)_{\lambda}=\{v\in {\text{Ind}}_{L(\frak p)}^{L(\frak g)}V_{\frak k}(\lambda): hv=\lambda(h)v\ \forall \ h\in\frak h\}\end{equation*} has dimension one 
and hence ${\text{Ind}}_{L(\frak p)}^{L(\frak g)}V_{\frak k}(\lambda)$ has a unique irreducible quotient obtained as follows.   
Let $ev:\loopg\to\frak g$ be the Lie algebra homomorphism obtained by mapping $t\mapsto 1$.  This gives a $\loopg$-module structure on any $\frak g$-module, in particular on the irreducible $\frak g$-module $V(\lambda)$  with highest weight $\lambda$. 
It is easy to see that the $\loopg$-module $V(\lambda)$ is the unique irreducible quotient of ${\text{Ind}}_{L(\frak p)}^{L(\frak g)}V_{\frak k}(\lambda)$. 

Assume from now on that $\lambda(h_{\alpha_n}) =0$ (and $\lambda(h_{\alpha_{n-1}}) =0$ in the case of $so(2n)$).  It is easy to see in this case that the element $x^-_{\alpha_n}\otimes 1.(1\otimes v_\lambda)$ (and $x^-_{\alpha_{n-1}}\otimes 1.(1\otimes v_\lambda)$, in the case of $so(2n)$) generates a proper submodule of ${\text{Ind}}_{L(\frak p)}^{L(\frak g)}V_{\frak k}(\lambda)$. Denote by  $W^{\text{aff}}(\lambda)$  the corresponding $\loopg$-module quotient, and let $w_\lambda\in W^{\text{aff}}(\lambda)$ be the image of $1\otimes v_\lambda$. 

\begin{prop}\label{genrel} \hfill
\begin{enumerate}
\item[(i)] As an $L(\frak g)$-module, 
$W^{\text{\rm aff}}(\lambda)$ is  generated by $w_\lambda$ with the following relations:
 \begin{align*} x_\alpha^+\otimes t^r. w_\lambda=0,\ \ & \ \ h\otimes t^r.w_\lambda=\lambda(h)w_\lambda,\ \ \alpha\in R^+,\ h\in\frak h,\\ (x_{\alpha_i}^-)^{\lambda(h_i)+1}. w_\lambda =0,\ \ &\ \ x_\beta^-\otimes (t-1)t^r. w_\lambda =0,\ \forall\ \beta\in R^+(J),\ \ r\in\bz.\end{align*}
\item[(ii)] For all $\lambda\in P^+$, $\text{dim}\ W^{\text{\rm aff}}(\lambda)<\infty$. 
\end{enumerate}
\end{prop}
\begin{pf} Part (i) is clear from the definition of $W^{\text{aff}}(\lambda)$. Part (ii) was proved in  \cite{CPweyl}. In the language of that paper, it is easy to see that the module $W^{\text{aff}}(\lambda)$ is a  quotient of $W(\bpi)$, where $\bpi=(\pi_1,\pi_2,\cdots ,\pi_n)$ is such that  $\pi_i(u)=(1-u)^{\lambda(h_{\alpha_i})}$ for all $1\le i\le n$.
\end{pf}

In view of the preceding proposition, we can write
\begin{equation*} W^{\text{aff}}(\lambda)=\bigoplus_{\mu\in P^+}m_{\lambda ,\mu} V(\mu).\end{equation*} Clearly, $m_{\lambda,\lambda} =1$ and $m_{\lambda,\mu}= 0$ if $\lambda-\mu\notin Q^+$.
We can now state our conjecture.

\begin{conj}\label{kc}
As $\frak g$-modules, we have
\begin{equation*} W^{\text{aff}}(\lambda)\cong W(\lambda).\end{equation*}
\end{conj}
Here $W(\lambda)$ is the representation $W_O(\lambda)$ or
$W_{Sp}(\lambda)$ of $\frak g$ defined in Section~\ref{ss_defW}.

\subsection{Type A sublattices}
As a first step towards providing evidence for this conjecture we prove the following result which is analogous to Proposition \ref{motive}.

For $\eta=\sum_i s_i\alpha_i\in Q^+$,  set $\text{supp}\ \eta=\{\alpha_i: s_i\ne 0\}$ 
and let $J(\eta)$ be the minimal connected subset of the set of simple roots that contains $\text{supp}\ \eta$. Let $\frak g(J_\eta)$ be the subalgebra generated by the elements $x_{\pm\alpha}$, $\alpha\in R^+(J_\eta)$. We say that $\eta$ is of type $A$ if the elements $x^\pm_{\alpha_i}$, $\alpha_i\in J(\eta)$ generate a subalgebra of $\frak g$ of type $sl_r$.
\begin{prop}  For all $0\ne \eta\in Q^+$ of type $A$, we have
\begin{equation*} m_{\lambda, \lambda-\eta} =0.\end{equation*}
\end{prop}

\begin{pf} Suppose first that $\alpha_n\notin\text{supp} \ \eta$.
Then $\text{supp}\ \eta\subset J$ and hence \begin{equation*} W^{\text{aff}}(\lambda)_{\lambda-\eta}= 1\otimes V_{\frak k}(\lambda)_{\lambda-\eta} .\end{equation*}
Since $V_{\frak k} (\lambda)$ is an irreducible $\frak k$-module, it follows that $w=1\otimes v_\lambda$.

If $\alpha_n\in\text{supp}\ \eta$, then since $\eta$ is of type $A$, one can conclude by a simple inspection that one of the following must hold.
\vskip 6pt 

\noindent {\bf{Case 1.}} \ If $\frak g$ is of type $sp(2n)$ or $so(2n+1)$, then $\text{supp}\ \eta =\{\alpha_n\}$ and hence $\eta= s\alpha_n$ for some $s\ge 0$. By definition, we have
\begin{equation*} x_{\alpha_n}^-\otimes 1.w_\lambda =0.\end{equation*}
Applying $h_{\alpha_n}\otimes t^r$ to the above equation, we see again from the definition that
\begin{equation*} x_{\alpha_n}^-\otimes t^r. w_\lambda =0.\end{equation*}
Hence $W^{\text{aff}}(\lambda)_{\lambda-\eta} =0$ thus proving the proposition.

\noindent {\bf{Case 2. }} \ Suppose that $\frak g$ is of type $so(2n)$. 
If $\alpha\in R^+(\text{supp}\ \eta)\backslash R^+(J)$, then  it is easy to see that  there exists $\beta\in R^+(J_\eta)$ such that $\alpha=\beta+\alpha_n$ or $\alpha=\beta+\alpha_{n-1}$ or $\alpha=\beta+\alpha_{n-1}+\alpha_n$. Writing $x^-_{\alpha} =[x^-_\beta, x_{\alpha_n}^-]$ etc. we see in all cases that 
\begin{equation*} x_{\alpha}^-\otimes (t-1)t^r. w_\lambda =0,\ \ \forall r\in\bz. \end{equation*}
Hence $\bu(L(\frak g(J_\eta))).w_\lambda =\bu(\frak g(J_\eta)).w_\lambda$. Since $\bu(\frak g(J_\eta)).w_\lambda$ is an irreducible $\frak g(J_\eta)$-module, and  $W^{\text{aff}}(\lambda)_{\lambda-\eta}\subset \bu(\frak g(J_\eta)).w_\lambda$, the proposition follows.
\end{pf}

\subsection{Roots to consider}
In this section we examine the structure of $\loopg$ to determine
which $V(\mu)$ can possibly appear as summands in
$W^{\text{aff}}(\lambda)$, or equivalently, which $m_{\lambda,\mu}$
might be nonzero.  The main result is Proposition~\ref{nec}, which
says that $\lambda-\mu$ must lie in the $\bz_+$-span of a special
subset of the positive roots.

 Although the module $W^{\text{aff}}(\lambda)$ is not an evaluation
module,  the next result shows that it is in fact a module for the quotient of $\loopg$ by the ideal $\frak g\otimes (t-1)^2\bc[t,t^{-1}]$.

\begin{prop}{\label{evmod}}\hfill

\begin{enumerate}
\item[(i)] Let $\alpha\in R^+$ be such that $\alpha-\alpha_n\in R^+(J)$ (or $\alpha-\alpha_{n-1}\in R^+(\overline{J})$ when $\frak g=so(2n)$). 
Then for all $r\in\bz$ with $|r|\ge 1$,
\begin{equation*} (x_\alpha^-\otimes (t-1)^r). w_\lambda=0.\end{equation*}
\item[(ii)]  For all $\alpha\in R^+$, and $r\in\bz$ with $|r|\ge 2$, we have,
\begin{equation*} (x_\alpha^-\otimes (t-1)^r). w_\lambda=0.\end{equation*}
\end{enumerate}
\end{prop}
\begin{pf} To prove (i), recall that
\begin{equation*}x_{\alpha}^-\otimes (t-1)^r =[x_{\alpha_n}^-\otimes 1,  x_{\alpha-\alpha_n}^-\otimes (t-1)^r].\end{equation*}
Part (i) now follows from the defining relations in $W^{\text{aff}}(\lambda)$.

It suffices to prove (ii) in the case of $\theta$, the highest root.
A simple case by case inspection shows that we can write  $\theta =\alpha+\beta$ where either
\begin{enumerate}
\item[(a)]  $\alpha\in R^+(J)$,  $\beta\in R^+$, $\beta-\alpha_n\in R^+(J)$, or \item[(b)] both $\alpha,\beta\in R^+$ and $\alpha-\alpha_n,\ \beta-\alpha_n\in R^+(J)$,
\item[(c)] both $\alpha,\beta\in R^+$ and $\alpha-\alpha_n,\in R^+(J)$, $\beta-\alpha_{n-1}\in R^+(\overline{J})$.
\end{enumerate}
The result now follows from part (i).
\end{pf}

In view of the preceding proposition, we shall be interested in the following subset of $R^+$,
\begin{equation*}\{\beta\in R^+: \beta, \beta-\alpha_n\notin R^+ \setminus R^+(J)\}.\end{equation*}
(or, if $\frak g=so(2n)$, the set \begin{equation*}\{\beta\in R^+:
\beta, \beta-\alpha_n,\beta-\alpha_{n-1}\notin R^+ \setminus \left(R^+(J)\cup R^+(\overline {J})\right\}).\end{equation*} Let $N$ be the cardinality of this set. 

We write the sets out explicitly for the reader's convenience. The sets in questions can be written as the collection of roots   $\beta_{k,l}$, defined as follows. 
\vskip 6pt
\noindent $\frak g=sp(2n)$. For $1\le k<l\le n-1$, set \begin{equation*} \beta_{k,l} = \alpha_k+\alpha_{k+1}+\cdots +\alpha_{l-1}+2\alpha_{l}+\cdots +2\alpha_{n-1}+\alpha_n,\end{equation*}
and set $\beta_{0,l} =2\alpha_1+\cdots+ 2\alpha_{n_1}+\alpha_n$.

\vskip 6pt
 \noindent $\frak g= so(2n+1)$. For $1\le k<l\le n-1$, set \begin{equation*} \beta_{k,l}= \alpha_k+\alpha_{k+1}+\cdots +\alpha_{l-1}+2\alpha_l+\cdots +2\alpha_n.\end{equation*}

\vskip 6pt
\noindent 
$\frak g=so(2n+2)$. For $1\le k<l\le n-1$, set
\begin{equation*} \beta_{k,l}= \alpha_k+\alpha_{k+1}+\cdots +\alpha_{l-1}+2\alpha_l+\cdots +2\alpha_{n-1}+\alpha_{n}+\alpha_{n+1}.\end{equation*}

We order the roots $\beta_{k,l}$ in the lexicographic order,\begin{equation*}\beta_{k,l} < \beta_{k',l'}\ \text{if} \ k<k',\ \text{or}\  k=k', l<l',\end{equation*} and let $\beta_1,\beta_2,\cdots ,\beta_N$ be an enumeration of these roots. 

Given $\bor\in\bz^N$, set
\begin{equation*} x^-_\bor =(x^-_{\beta_1}\otimes (t-1))^{r_1}(x^-_{\beta_2}\otimes (t-1))^{r_2}\cdots 
(x^-_{\beta_N}\otimes (t-1))^{r_N}.\end{equation*}
Let $\le $ be the lexicographic order on $\bz^N$.  The following corollary to Proposition \ref{evmod} is now immediate, by a simple application of the PBW theorem.

\begin{cor}\label{span} We have,
\begin{equation*} W^{\text{\rm aff}}(\lambda) =\sum_{\bor\in\bz_+^{N}}\bu(\frak n^-)x_\bor^-.w_\lambda\end{equation*}
\hfill\qedsymbol\end{cor}

We can now give a necessary condition for $m_{\lambda,\mu}$ to be nonzero. We begin with the following Lemma which is 
 easily checked.
\begin{lem} \label{commute}
\begin{enumerate} \item[(i)] For all $1\le r,s\le N$ we have $\beta_r+\beta_s\notin R^+$.
\item[(ii)] For any $i=1,\ldots,n$ and for all $1\le r, s\le N$, we have
$\beta_r+\beta_s-\alpha_i\notin R^+$.
\item[(iii)] Suppose that $i=1,\ldots,n$ and $1\le r\le N$ are such that $\beta_r-\alpha_i\in R^+$. Then, either $\beta_r-\alpha_i =\beta_s$ for some $s\ge r$, or $\beta_r=\beta_{k,n-1}$, $i=n-1$ for some $k$.\end{enumerate}
\hfill\qedsymbol\end{lem}

\begin{prop}\label{nec} Assume that $m_{\lambda,\mu}\ne 0$. Then
\begin{equation*}\mu=\lambda-\sum_{j=1}^N s_j\beta_j,\end{equation*}
for some nonnegative integers $s_1,s_2,\cdots s_N$.
\end{prop}

\begin{pf} Let $W_1$ be a $\frak g$-module complement to $V(\lambda)$ so that we have
\begin{equation*} W^{\text{aff}}(\lambda)= V(\lambda)\oplus W_1,\end{equation*}  as $\frak g$-modules. If $W_1\ne 0$, choose $\bor_1$ minimal so that the projection $w_{\bor_1}$ of $x^-_{\bor_1}.w_\lambda$ onto $W_1$ is nonzero. Using Lemma \ref{commute}, we see that  
\begin{equation*} x_{\alpha_i}^+.  x_{\bor_1}.w_\lambda \in\sum_{\bos<\bor_1}\bu(\frak n^-).x^-_\bos. w_\lambda. \end{equation*} The minimality of $\bor_1$ now implies that 
\begin{equation*} x_{\alpha_i}^+.  x_{\bor_1}.w_\lambda\in V(\lambda),\end{equation*} which implies that \begin{equation*} x_{\alpha_i}^+.w_{\bor_1}.w_\lambda=0.\end{equation*}
Hence, $w_{\bor_1}$  generates an irreducible $\frak g$-module. Let $W_2\subset W_1$ be the $\frak g$-module complement to it. Repeating the argument, we see that there exist a finite set $\bor_1,\bor_2,\cdots \bor_m$ such that 
\begin{equation*} W^{\text{aff}}(\lambda) = \oplus V(\mu_{\bor_j}),\end{equation*}
where $\mu_{\bor_j}$ is the weight of the element $w_{\bor_j}$. Clearly, each $\mu_{\bor_j}$ has weight of the form $\lambda-\sum_j s_j\beta_j$ and the proposition follows.
\end{pf}

\subsection{Stable limit property}

We next prove an interesting consequence of the preceding
proposition. It is the analogue of the statement following
Definition~\ref{def_W} that the modules $W(\lambda)$ have the same
decomposition for all sufficiently large $n$. We prove this for the
orthogonal algebras; the case of symplectic algebras is similar and
simpler.

For this proposition only, we denote by $\frak {g}_{2n}$ the Lie algebra   $so(2n)$ and by $\frak{g}_{2n+1}$ the Lie algebra $so(2n+1)$. We denote the corresponding lattice $Q^+$ by $Q_n^+$ etc., and similarly denote $W^{\text{aff}}(\lambda)$ by $W_n^{\text{aff}}(\lambda) $,  and the multiplicities $m_{\lambda,\mu}$ by $m_{\lambda,\mu, n}$.

We  have an embedding of $so(2n-1)\to so(2n)$, given as follows,
\begin{equation*} x_{\alpha_i}^\pm\mapsto x^{\pm}_{\alpha_i}, \ \ 1\le i\le n-2,\ \ x^\pm_{\alpha_{n-1}}\mapsto x_{\alpha_n}+x_{\alpha_{n-1}}.\end{equation*}
In other words $so(2n-1)$ is the subalgebra of fixed points of the automorphisms of $so(2n)$  defined by interchanging the spin nodes of the Dynkin diagram..  
Under this embedding, the root vectors
\begin{equation*} x^-_{\beta_{k,l}}\mapsto  x^\pm_{\beta_{k,l}}, \ \ 1\le k<l\le n-1.\end{equation*}  
The restriction map $\frak{h}^*_{2n}\to\frak{h}^*_{2n-1}$ induces an isomorphism between the subspace of $P_{2n}^+$  spanned by $\omega_i$, $1\le i\le n-2$ and the subspace of $P^+_{2n-1}$ spanned by $\omega_i$, $1\le i\le n-2$.

We also define an embedding of $so(2n)\to so(2n+1)$. This is given by the assignment,
\begin{align*} x^\pm_{\alpha_i}\mapsto x^\pm_{\alpha_i}, &\ \ 1\le i\le n-2,\\
x^\pm_{\alpha_{n-1}}\mapsto x^\pm_{\alpha_{n-1}}, \ \ &\ \ 
x^\pm_{\alpha_{n}}\mapsto x^\pm_{\alpha_{n-1}+2\alpha_{n}}.\end{align*}
Again notice that this embedding maps $\beta_{k,l}$ to $\beta_{k,l}$
for all $l\le n-2$ and $\omega_i$ to $\omega_i$ for $i\le n-2$. 

Both  embeddings naturally extend to maps of the corresponding loop algebras.

\begin{thm} Let $\lambda\in P_n^+$ and assume that $\lambda(h_n) =\lambda(h_{n-1})=0$. There exists $r(n)\ge n$ such that for all $s, s'\ge r(n)$, we have
\begin{equation*} m_{\lambda,\mu,s}=m_{\lambda,\mu, s'}.\end{equation*}
\end{thm}
The theorem is clearly a consequence of the following proposition.
\begin{prop} Let $\lambda\in P_n^+$. Then,
\begin{equation*} m_{\lambda,\mu,n}\ge m_{\lambda, \mu,n+1}.\end{equation*}
\end{prop}
\begin{pf} Assume first that $n=2m+1$ and let $N$ be the number of  roots of the form $\beta_{i,k}$ for $so(2m+1)$. Notice that this is exactly the same number of such roots for $so(2m+2)$. By Proposition \ref{span}, we have
\begin{equation*} W_{2m+2}^{\text{aff}}(\lambda) =\sum_{\bor\in\bz^N}\bu(\frak{n}^-_{2m+2})x^-_\bor.w_\lambda.\end{equation*}
The elements $x^-_{\bor}$ are in the image of the embedding of $so(2m+1)\to so(2m+2)$. Hence,  setting, $W_{2m+1} =\bu(\frak{g}_{2m+1}).w_\lambda$ we see that
\begin{equation*} W_{2m+1} =\sum_{\bor\in\bz^N} \bu(\frak{n}^-_{2m+1}).w_\lambda.\end{equation*}
It is now clear, that the elements $w_\bos$ defined in the proof of Proposition \ref{nec}  can be chosen to be in $W_{2m+1}$. It is clear from the defining relations of $W^{\text{aff}}_{2m+1}(\lambda)$ that there exists a surjective map  $W^{\text{aff}}_{2m+1}(\lambda)\to W_{2m+1}$ of $\frak{g}_{2m+1}$--modules.
This clearly implies that $m_{\mu, 2m+2}\le m_{\mu, 2m+1}$.

To prove that $m_{\mu, 2m+1}\le m_{\mu, 2m}$, we  use the embedding of $so(2m)\to so(2m+1)$. The proof is similar, we just need to show that the elements $x^-_\bor$ that span $W_{2m+1}^{\text{aff}}(\lambda)$ are actually in $\bu(\frak{n}_{2m})$. The only difficulty is with the roots $\beta_{k,n-1}$ and $\beta_{k, n}$. Now, $\beta = \beta_{k,2n-1}-\alpha_{2m-1}-\alpha_n\in R^+$ is such that 
$x_\beta^-\otimes (t-1).w_\lambda=0$. Since $\lambda(h_{n-1})=\lambda(h_n) =0$, it follows that $x^-_{\alpha_{n-1}+\alpha_n}=0$. This now implies that 
$x_{\beta_{k,n-1}}^-\otimes (t-1).w_\lambda=0$. The case of $\beta_{k,n}$ is dealt with similarly.
This completes the proof of the proposition.

\end{pf}

}
{\newcommand{\thmref}[1]{Theorem~\ref{#1}}
\newcommand{\secref}[1]{Sect.~\ref{#1}}
\newcommand{\lemref}[1]{Lemma~\ref{#1}}
\newcommand{\propref}[1]{Proposition~\ref{#1}}
\newcommand{\corref}[1]{Corollary~\ref{#1}}
\newcommand{\remref}[1]{Remark~\ref{#1}}
\newcommand{\nc}{\newcommand}
\newcommand{\rnc}{\renewcommand}
\nc{\cal}{\mathcal}
\nc{\goth}{\mathfrak}
\rnc{\bold}{\mathbf}
\renewcommand{\frak}{\mathfrak}
\renewcommand{\Bbb}{\mathbb}

\nc{\Cal}{\mathcal}
\nc{\Xp}[1]{X^+(#1)}
\nc{\Xm}[1]{X^-(#1)}
\nc{\on}{\operatorname}
\nc{\ch}{\mbox{ch}}
\nc{\Z}{{\bold Z}}
\nc{\J}{{\cal J}}
\nc{\C}{{\bold C}}
\nc{\Q}{{\bold Q}}
\renewcommand{\P}{{\cal P}}
\nc{\N}{{\Bbb N}}
\nc\beq{\begin{equation}}
\nc\enq{\end{equation}}
\nc\lan{\langle}
\nc\ran{\rangle}
\nc\bsl{\backslash}
\nc\mto{\mapsto}
\nc\lra{\leftrightarrow}
\nc\hra{\hookrightarrow}
\nc\sm{\smallmatrix}
\nc\esm{\endsmallmatrix}
\nc\sub{\subset}
\nc\ti{\tilde}
\nc\nl{\newline}
\nc\fra{\frac}
\nc\und{\underline}
\nc\ov{\overline}
\nc\ot{\otimes}
\nc\bbq{\bar{\bq}_l}
\nc\bcc{\thickfracwithdelims[]\thickness0}
\nc\ad{\text{\rm ad}}
\nc\Ad{\text{\rm Ad}}
\nc\Hom{\text{\rm Hom}}
\nc\End{\text{\rm End}}
\nc\Ind{\text{\rm Ind}}
\nc\Res{\text{\rm Res}}
\nc\Ker{\text{\rm Ker}}
\rnc\Im{\text{Im}}
\nc\sgn{\text{\rm sgn}}
\nc\tr{\text{\rm tr}}
\nc\Tr{\text{\rm Tr}}
\nc\supp{\text{\rm supp}}
\nc\card{\text{\rm card}}
\nc\bst{{}^\bigstar\!}
\nc\he{\heartsuit}
\nc\clu{\clubsuit}
\nc\spa{\spadesuit}
\nc\di{\diamond}

\nc\al{\alpha}
\nc\bet{\beta}
\nc\ga{\gamma}
\nc\de{\delta}
\nc\ep{\epsilon}
\nc\io{\iota}
\nc\om{\omega}
\nc\si{\sigma}
\rnc\th{\theta}
\nc\ka{\kappa}
\nc\la{\lambda}
\nc\ze{\zeta}

\nc\vp{\varpi}
\nc\vt{\vartheta}
\nc\vr{\varrho}

\nc\Ga{\Gamma}
\nc\De{\Delta}
\nc\Om{\Omega}
\nc\Si{\Sigma}
\nc\Th{\Theta}
\nc\La{\Lambda}
\nc\boa{\bold a}
\nc\bob{\bold b}
\nc\boc{\bold c}
\nc\bod{\bold d}
\nc\boe{\bold e}
\nc\bof{\bold f}
\nc\bog{\bold g}
\nc\boh{\bold h}
\nc\boi{\bold i}
\nc\boj{\bold j}
\nc\bok{\bold k}
\nc\bol{\bold l}
\nc\bom{\bold m}
\nc\bon{\bold n}
\nc\boo{\bold o}
\nc\bop{\bold p}
\nc\boq{\bold q}
\nc\bor{\bold r}
\nc\bos{\bold s}
\nc\bou{\bold u}
\nc\bov{\bold v}
\nc\bow{\bold w}
\nc\boz{\bold z}

\nc\ba{\bold A}
\nc\bb{\bold B}
\nc\bc{\bold C}
\nc\bd{\bold D}
\nc\be{\bold E}
\nc\bg{\bold G}
\nc\bh{\bold h}
\nc\bH{\bold H}

\nc\bi{\bold I}
\nc\bj{\bold J}
\nc\bk{\bold K}
\nc\bl{\bold L}
\nc\bm{\bold M}
\nc\bn{\bold N}
\nc\bo{\bold O}
\nc\bp{\bold P}
\nc\bq{\bold Q}
\nc\br{\bold R}
\nc\bs{\bold S}
\nc\bt{\bold T}
\nc\bu{\bold U}
\nc\bv{\bold v}
\nc\bV{\bold V}

\nc\bw{\bold W}
\nc\bz{\bold Z}
\nc\bx{\bold x}
\nc\bX{\bold X}
\nc\blambda{{\mbox{\boldmath $\Lambda$}}}
\nc\bpi{{\mbox{\boldmath $\pi$}}}

\newcommand{\vdomsize}[2]{{\begin{picture}(#1,#2)\multiput(0,0)(0,#1){3}%
{\line(1,0){#1}}\multiput(0,0)(#1,0){2}{\line(0,1){#2}}\end{picture}}}
\newcommand{\hdomsize}[2]{{\begin{picture}(#2,#1)\multiput(0,0)(#1,0){3}%
{\line(0,1){#1}}\multiput(0,0)(0,#1){2}{\line(1,0){#2}}\end{picture}}}
\newcommand{\vdom}{{\vdomsize48}} \newcommand{\vdoms}{{\vdomsize36}}
\newcommand{\hdom}{{\hdomsize48}} \newcommand{\hdoms}{{\hdomsize36}}

\nc\e[1]{E_{#1}}
\nc\ei[1]{E_{\delta - \alpha_{#1}}}
\nc\esi[1]{E_{s \delta - \alpha_{#1}}}
\nc\eri[1]{E_{r \delta - \alpha_{#1}}}
\nc\ed[2][]{E_{#1 \delta,#2}}
\nc\ekd[1]{E_{k \delta,#1}}
\nc\emd[1]{E_{m \delta,#1}}
\nc\erd[1]{E_{r \delta,#1}}

\nc\ef[1]{F_{#1}}
\nc\efi[1]{F_{\delta - \alpha_{#1}}}
\nc\efsi[1]{F_{s \delta - \alpha_{#1}}}
\nc\efri[1]{F_{r \delta - \alpha_{#1}}}
\nc\efd[2][]{F_{#1 \delta,#2}}
\nc\efkd[1]{F_{k \delta,#1}}
\nc\efmd[1]{F_{m \delta,#1}}
\nc\efrd[1]{F_{r \delta,#1}}
\nc{\ug}{\bu^{fin}}

\nc\fa{\frak a}
\nc\fb{\frak b}
\nc\fc{\frak c}
\nc\fd{\frak d}
\nc\fe{\frak e}
\nc\ff{\frak f}
\nc\fg{\frak g}
\nc\fh{\frak h}
\nc\fj{\frak j}
\nc\fk{\frak k}
\nc\fl{\frak l}
\nc\fm{\frak m}
\nc\fn{\frak n}
\nc\fo{\frak o}
\nc\fp{\frak p}
\nc\fq{\frak q}
\nc\fr{\frak r}
\nc\fs{\frak s}
\nc\ft{\frak t}
\nc\fu{\frak u}
\nc\fv{\frak v}
\nc\fz{\frak z}
\nc\fx{\frak x}
\nc\fy{\frak y}

\nc\fA{\frak A}
\nc\fB{\frak B}
\nc\fC{\frak C}
\nc\fD{\frak D}
\nc\fE{\frak E}
\nc\fF{\frak F}
\nc\fG{\frak G}
\nc\fH{\frak H}
\nc\fJ{\frak J}
\nc\fK{\frak K}
\nc\fL{\frak L}
\nc\fM{\frak M}
\nc\fN{\frak N}
\nc\fO{\frak O}
\nc\fP{\frak P}
\nc\fQ{\frak Q}
\nc\fR{\frak R}
\nc\fS{\frak S}
\nc\fT{\frak T}
\nc\fU{\frak U}
\nc\fV{\frak V}
\nc\fZ{\frak Z}
\nc\fX{\frak X}
\nc\fY{\frak Y}
\nc\tfi{\ti{\Phi}}
\nc\bF{\bold F}

\nc\ua{\bold U_\A}

\nc\qinti[1]{[#1]_i}
\nc\q[1]{[#1]_q}
\nc\xpm[2]{E_{#2 \delta \pm \alpha_#1}}  
\nc\xmp[2]{E_{#2 \delta \mp \alpha_#1}}
\nc\xp[2]{E_{#2 \delta + \alpha_{#1}}}
\nc\xm[2]{E_{#2 \delta - \alpha_{#1}}}
\nc\hik{\ed{k}{i}}
\nc\hjl{\ed{l}{j}}
\nc\qcoeff[3]{\left[ \begin{smallmatrix} {#1}& \\ {#2}& \end{smallmatrix}
\negthickspace \right]_{#3}}
\nc\qi{q}
\nc\qj{q}

\nc\ufdm{{_\ca\bu}_{\rm fd}^{\le 0}}


\nc\isom{\cong} 

\nc{\pone}{{\Bbb C}{\Bbb P}^1}
\nc{\pa}{\partial}
\def\H{\cal H}
\def\L{\cal L}
\nc{\F}{{\cal F}}
\nc{\Sym}{{\goth S}}
\nc{\A}{{\cal A}}
\nc{\arr}{\rightarrow}
\nc{\larr}{\longrightarrow}

\nc{\ri}{\rangle}
\nc{\lef}{\langle}
\nc{\W}{{\cal W}}
\nc{\uqatwoatone}{{U_{q,1}}(\su)}
\nc{\uqtwo}{U_q(\goth{sl}_2)}
\nc{\dij}{\delta_{ij}}
\nc{\divei}{E_{\alpha_i}^{(n)}}
\nc{\divfi}{F_{\alpha_i}^{(n)}}
\nc{\Lzero}{\Lambda_0}
\nc{\Lone}{\Lambda_1}
\nc{\ve}{\varepsilon}
\nc{\phioneminusi}{\Phi^{(1-i,i)}}
\nc{\phioneminusistar}{\Phi^{* (1-i,i)}}
\nc{\phii}{\Phi^{(i,1-i)}}
\nc{\Li}{\Lambda_i}
\nc{\Loneminusi}{\Lambda_{1-i}}
\nc{\vtimesz}{v_\ve \otimes z^m}

\nc{\asltwo}{\widehat{\goth{sl}_2}}
\nc\eh{\frak h^e}  
\nc\loopg{L(\frak g)}  
\nc\eloopg{L^e(\frak g)} 
\nc\ebu{\bu^e} 
\nc\loopa{L(\frak a)}  

\nc\teb{\tilde E_\boc}
\nc\tebp{\tilde E_{\boc'}}

\section{Motivations and special cases of the conjecture}  We have restricted ourselves to the case of the enveloping algebras of loop algebras to simplify matters and to avoid excessive notation. However, the motivation for Conjecture \ref{kc} comes from connections which we now explain, with the irreducible finite-dimensional   representations of quantum affine algebras. 

\subsection{Background}
\label{sec_kr}
\newcommand{\g}{{\mathfrak g}}
\renewcommand{\l}{\ell}

The quantum affine algebra $U_q(\hat{\g})$ and the related Yangian $Y(\g)$
were introduced by Drinfeld and Jimbo as tools for studying solutions
to the quantum Yang--Baxter equation.  Finite dimensional
representations of either Hopf algebra give rise to solutions (called
$R$-matrices), which can, in turn, be used to construct the transfer
matrices of integrable dynamical systems.  The Bethe Ansatz is a
technique for calculating eigenvalues of such transfer matrices, as
the solutions to a set of algebraic equations.

The algebras above have subalgebras $U_q(\g)\hookrightarrow
U_q(\hat{\g})$ and $\g\hookrightarrow Y(\g)$, and the eigenspaces of
the transfer matrix decompose a finite-dimensional representation of
the larger algebra into subspaces stabilized by the smaller one.  If
each eigenspace were a single irreducible representation, then the
eigenvalues of the transfer matrix would completely describe the
decomposition of a representation of $U_q(\hat{\g})$ or $Y(\g)$ upon
restriction to $U_q(\g)$ or $\g$, respectively.  If, in addition, the
Bethe Ansatz finds all eigenvalues, then solving the Bethe equations
would yield the complete desired decomposition.

This was the approach employed by Kirillov and Reshetikhin~\cite{KR}.
They addressed the problem of decomposing an irreducible $Y(\g)$
module according to the action of the embedded copy of $\g$ by
conjecturing that the Bethe Ansatz detected all the pieces in the
decomposition.  The result was a so-called ``fermionic formula'' for
the number of times each $\g$ module would appear in the
decomposition.  Their attention was restricted to a particular class
of finite-dimensional representations in which the Bethe eigenvectors
are especially well-behaved.

The results in this paper are all from the point of view of the
embedding $U_q(\g)\hookrightarrow U_q(\hat{\g})$.  It has long been a
folk theorem that the decompositions in this case were identical to
those in the $\g\hookrightarrow Y(\g)$ case.  A proof was recently
given for simply-laced $\g$ by Varagnolo~\cite{V}.  Further, we know
by results of Lusztig~\cite{L} that the representation theory of
$\bu_q(\frak g)$ over $\bc(q)$ is the `same' as the
representation theory of $\bu(\frak g)$.  Hence we we are justified in
talking about the Kirillov--Reshetikhin conjecture on Yangians as if
it applies to the representations we studied in Section~\ref{sec_algebra}.

\vskip 6pt

{\bf{Conjecture of Kirillov and Reshetikhin.}}
For each $m\in\bz^+$ and $\l=1,\ldots,n$, there exists  an irreducible representation $V_q(m\omega_\l)$  of $\bu_q(\hat{\g})$ whose
highest weight when viewed as a representation of $\bu_q(\g)$ is
$m\omega_\l$.  Further, the decomposition of the tensor product of $N$ such 
representations as $\bu_q(\frak g)$-modules is given by
$$
\bigotimes_{a=1}^N (V_q(m_a\omega_{\l_a}) |_{\g})
\simeq \sum_{\lambda} n_\lambda V(\lambda)
$$
where the sum runs over all weights $\lambda$ less than $\sum m_a\omega_{\l_a}$,
the highest weight of the tensor product.  The nonnegative integer
$n_\lambda$ is the multiplicity with which the irreducible $\g$-module
$V(\lambda)$ occurs.  Write $\lambda =  \sum m_a\omega_{\l_a} - \sum
n_i \alpha_i$.  Then
$$
n_\lambda = 
\sum_{\mbox{partitions}} \;\; \prod_{n\geq1} \;\; \prod_{k=1}^r \;\;
\left(\!\!\begin{array}{cc}P^{(k)}_n(\nu) + \nu^{(k)}_n \\ \nu^{(k)}_n
       \end{array}\!\!\right)
$$
The sum is taken over all ways of choosing partitions
$\nu^{(1)},\ldots,\nu^{(r)}$ such that $\nu^{(i)}$ is a partition of
$n_i$ which has $\nu^{(i)}_n$ parts of size $n$ (so $n_i =
\sum_{n\geq1} n \nu^{(i)}_n$).  The function $P$ is defined by
\begin{eqnarray*}
P^{(k)}_n(\nu) &=& \sum_{a=1}^N \min(n,m_a)\delta_{k,\l_a}
  - 2 \sum_{h\geq 1} \min(n,h)\nu^{(k)}_{h} + \\
&&\hspace{1cm} +
  \sum_{j\neq k}^r \sum_{h\geq 1} \min(-c_{k,j}n,-c_{j,k}h)\nu^{(j)}_{h}
\nonumber
\end{eqnarray*}
where $C=(c_{i,j})$ is the Cartan matrix of $\g$, and $\binom{a}{b}=0$
whenever $a<b$.

\vskip 6pt 
The formula describing the $n_\lambda$  is called the fermionic
formula.  The connection with representation theory was made by
Kirillov and Reshetikhin, who proved the conjecture in the case of $sl_n$.

The fermionic formula is somewhat difficult to work with directly.
However, when $\g$ is classical and we ignore the tensor product
(taking $N=1$), there is a simple combinatorial description of what
the fermionic formula predicts for $V_q(m\omega_\l)$; for a derivation
of the combinatorics from the fermionic formula see~\cite{kleber}.  In
particular, as long as the weight $\omega_\l$ lies in the type $A$
part of the Dynkin diagram ($\l<n$ and $\l<n-1$ for $so(2n)$), the
Kirillov--Reshetikhin decompositions are a special case of the ones we
defined in Section~\ref{sec_sym}: the $\bu_q(\frak g)$ module
structure of $V_q(m\omega_\l)$ is the same as the $\frak g$-module
structure of $W_O(m\omega_\l)$ when $\g$ is orthogonal, and is the
same as $W_{Sp}(m\omega_\l)$ when $\g$ is symplectic.  Thus it becomes
natural to make the following conjecture.

\begin{conj} \label{gkr}
 There exists an irreducible representation $V^{\text{aff}}_q(\lambda)$ of the quantum affine algebra whose $\bu_q(\frak g)$-module decomposition is $W_G(\lambda)$.  
\end{conj}

\subsection{Minimal affinizations} It is known~\cite{CP} that the
finite-dimensional irreducible representations of quantum affine
algebras $\bu_q(\frak g)$ are indexed by the $n$-tuples $(\pi_1,\ldots
,\pi_n)$ of polynomials with constant term 1. In \cite{C}, we  showed
that the module $V_q(m\omega_\l)$ conjectured by Kirillov and Reshetikhin is given by the $n$-tuple
\begin{equation*} \pi_j=1,\ \ j\ne m,\qquad \pi_m=(1-q^{-\l+1}u)(1-q^{-\l+3}u)\cdots (1-q^{\l-1}u).\end{equation*}
These modules are the so called  minimal affinization of $m\omega_\l$,
see \cite{C}. However, minimal affinizations
$V_q^{\text{aff}}(\lambda)$ are known to exist more generally for any
dominant integral weight $\lambda$ \cite{CP}, and as long as $\lambda$
is not supported on the spin nodes, as in the previous section, they
are  unique (up to $\bu_q(\frak g)$-module isomorphism).
 In fact, these  modules are determined by the requirement that
\begin{equation*} m_{\lambda,\lambda-\eta} =0 \quad \mbox{for all
$\eta\neq0$ not supported on the spin nodes,}\end{equation*}
a condition satisfied by the modules in Sections~\ref{sec_sym} and~\ref{sec_algebra} of this paper.

It can be shown as  in \cite{C} that 
 on specializing this representation, by putting $q=1$, we get a
quotient of the module $W^{\text{aff}}(\lambda)$. Thus, to prove our
Conjecture~\ref{gkr} generalizing the Kirillov--Reshetikhin decompositions, it suffices to prove Conjecture~\ref{kc} along with 
the statement that the specialized module is isomorphic to
$W^{\text{aff}}(\lambda)$. 

In the rest of this section, we restrict ourselves to the case of
$so(2n)$ and explicitly calculate $W(\lambda)=W_O(\lambda)$ for
several families of $\lambda$. We also show that the modules
$W^{\text{aff}}(\lambda)$ are isomorphic to a submodule of
$W(\lambda)$ in these cases. To prove that it is isomorphic to
$W(\lambda)$ requires arguments in the quantum algebra, similar to the
ones in \cite{C} and we do not give details of that here.

\subsection{Computation of examples}
\label{sec_compute}
Our computations use the technique described in Section~\ref{ss_combW}.
We begin with the simplest case.

\begin{example}[Rectangles]
\label{ex_Wrect}
When $\lambda$ is a multiple of a fundamental weight $m\omega_\ell$, its
Young diagram is a rectangle with $\ell$ rows and $m$ columns.  The
reader can readily verify that for each $\nu\in Y_\vdom$, there is a
unique semi-standard tableau for $\lambda/\nu$ whose reverse row
word is a ballot sequence:
$$
\setlength{\unitlength}{10pt}
\newcommand{\scribble}[2]{\multiput(#1,#2)(0,.33){6}{\line(3,1){1}}}
\newcommand{\num}[3]{\put(#1.5,#2.5){\makebox(0,0){\small$#3\cdots #3$}}}
\begin{picture}(9,5)
\put(0,0){\framebox(9,5){}}
\put(0,1){\framebox(3,2){}}
\put(0,3){\framebox(6,2){}}
\scribble01\scribble11\scribble21
\scribble03\scribble13\scribble23\scribble33\scribble43\scribble53
\num101\num421\num412\num403\num741\num732\num723\num714\num705
\end{picture}
$$
The Young diagram of its content $\mu$ is that of $\lambda/\nu$
rotated $180^\circ$.  Thus $W(m\omega_\ell)$ is the sum of $V(\mu)$ over 
all dominant weights $\mu$ which can be obtained from
$\lambda=m\omega_\ell$ by repeatedly subtracting some fundamental weight
$\omega_i$ and adding $\omega_{i-2}$ instead: each replacement of
$\omega_i$ by $\omega_{i-2}$ corresponds to a vertical domino in $\nu$
removing two boxes from a column of height $i$.

The proof that $W^{\text{aff}}(\lambda)\cong W_O(\lambda)$ was given in \cite{C}.

\end{example}

\begin{example}
\label{ex_Wabc}
Now take $\lambda=a\omega_1+b\omega_2+c\omega_3$; for $D_5$ this is
the generic weight not supported on the spin nodes.  This example is
simple because the only 
$\nu\in Y_\vdom$ that fit inside $\lambda$
consist of two rows of equal length.  The typical tableau this time is 
of the following form:
$$\setlength{\unitlength}{10pt}
\newcommand{\scribble}[2]{\multiput(#1,#2)(0,.33){6}{\line(3,1){1}}}
\newcommand{\num}[3]{\put(#1.5,#2.5){\makebox(0,0){\small$#3\cdots #3$}}}
\newcommand{\arrs}[3]{\put(#1.5,#2){\makebox(0,0){\small$#3$}}
\put(#1.1,#2){\vector(-1,0){1}}\put(#1.9,#2){\vector(1,0){1}}}
\begin{picture}(18,7)
\put(0,2){\line(0,1){3}}\multiput(12,2)(3,1){3}{\line(0,1){1}}
\put(0,2){\line(1,0){12}}\put(0,5){\line(1,0){18}}
\put(12,3){\line(1,0){3}}\put(15,4){\line(1,0){3}}
\put(0,3){\line(1,0){9}}\put(9,3){\line(0,1){2}}
\scribble03\scribble13\scribble23
\scribble33\scribble43\scribble53
\scribble63\scribble73\scribble83
\num121\num422\num723\num{10}23\num{10}32\num{10}41
\num{13}32\num{13}41\num{16}41
\multiput(0,.75)(3,0){4}{\line(0,1){.75}}
\put(0,5.5){\line(0,1){.75}}\multiput(12,5.5)(3,0){3}{\line(0,1){.75}}
\arrs11t\arrs41s\arrs71r\arrs{13}6b\arrs{16}6a
\put(6,6){\makebox(0,0){\small$c$}}
\put(5.6,6){\vector(-1,0){5.5}}\put(6.4,6){\vector(1,0){5.5}}
\end{picture}
$$

Here $\nu$ has two rows of length $r+s+t$, which is depicted here as
being less than $c$, but can also lie between $c$ and $b+c$.  In
either case, by definition
$$
s+t \leq c.
$$
Two additional restrictions on the parameters are imposed by the
requirement that the reverse row word be a ballot sequence:
\begin{align*}
s &\leq a \\
r &\leq b
\end{align*}
Requiring $r \leq b$ ensures that the numer of $2$s never exceeds the
number of $3$s; this is still the correct bound no matter how $r+s+t$
compares with $c$.

Converting the content of the tableau into a weight $\mu$, we find
that
$$
\mu=(a-s+t)\omega_1+(b-r+s)\omega_2+(c-s-t)\omega_3.
$$
Thus the decomposition for $\lambda=a\omega_1+b\omega_2+c\omega_3$ is
$$
W(\lambda) =
\sum_{\begin{array}{c} \scriptstyle s\leq a,\,\, r\leq b, \\
\scriptstyle s+t\leq c\end{array}}
V(\lambda - r(\omega_2) - s(\omega_3-\omega_2+\omega_1) - t(\omega_3-\omega_1))
$$
We have rewritten $\mu$ to highlight the fact that each of $r$, $s$
and $t$ count the number of times some weight is subtracted from
$\lambda$.

\end{example}

We now turn to the $\loopg$-modules $W^{\text{aff}}(\lambda)$ for $D_5$. In this case we take $J=\{\alpha_1,\alpha_2,\alpha_3,\alpha_4\}$. The roots $\beta_{k,l}$ of Section 2 can also be written as \begin{equation*}
\beta_{1,2}=\omega_2,\ \ 
\beta_{1,3} = \omega_1-\omega_2+\omega_3,\ \ 
\beta_{2,3}= \omega_3-\omega_1.\end{equation*}
Thus, Proposition  \ref{nec} gives us that
\begin{equation*} W^{\text{aff}}(\lambda) =\bigoplus m_{\lambda,\mu}V(\mu),\end{equation*}
where \begin{equation*} \mu\in\{\lambda-r_1\omega_2-r_2(\omega_3-\omega_2 +\omega_1)-r_3(\omega_3-\omega_1): r_1,r_2,r_3\ge 0\}.\end{equation*}
This immediately gives
\begin{equation*} m_{\lambda,\mu}\ne 0\implies c\ge r_2+r_3.\end{equation*}
 
It remains to prove that $r_1\le b$ and $r_2\le a$ if $m_{\lambda,\mu}\ne 0$. 
To do this it is obviously enough to establish the following lemma. 
\begin{lem}{\label{reduce}} Let $\bor\in \bz^3$. Then, \begin{equation*} x^-_\bor.w_\lambda\in \sum_{\bos\le (b,a, r_3)}\bu(\frak g).x^-_\bos.w_\lambda\end{equation*}
\end{lem}
\begin{pf} Set $N= \sum_{\bos\le (b,r_2, r_3)}\bu(\frak g).x^-_\bos.w_\lambda$. Observe that,
\begin{equation*}(x_{\alpha_2}^-)^a. x^-_\bor = x^-_\bor x_{\alpha_2}^- +x^-_{r_1+1, r_2-1, r-3+1}.\end{equation*} If $r_1+1\le b$ this implies that $x_\bor x_{\alpha_2}^-.w_\lambda \in N$.
Repeating this, we see that $x_\bor(x^-_{\alpha_2})^l\in N$ if $r_1+l\le b$. Taking $l=b$, $r_1 =0$, we get
$x_{0,r_2, r_3}(x^-_{\alpha_2})^b.w_\lambda\in N$. Since
$(x^-_{\alpha_2})^{b+1}.w_\lambda =0$, we can apply $x^-_{\alpha_2}$
to $x_{(0,r-2, r_3)}(x^-_{\alpha_2})^b.w_\lambda$ to find that
\begin{equation*} x_{(1, r_2-1, r_3)}.w_\lambda\in N\end{equation*}
for all $r_2$.  But $x_{(1, r_2, r_3)}(x_{\alpha_2}^-){a-1}.w_\lambda\in N$, hence we get $x_{(2, r_2-1, r_3)}.(x_{\alpha_2}^-)^{a-1}.w_\lambda\in N$. Continuing, we get
$x_{(r_1, r_2, r_3)}.w_\lambda\in N$.

A similar argument shows how to reduce to $r_2\le a$. We omit the details. 
\end{pf}

\begin{example}
\label{ex_W24}
All the decompositions $W(\lambda)$ calculated so far have been
multiplicity free.  This is not the case in general; for completeness
we include the minimal counterexample, $\lambda=\omega_2+\omega_4$.
This time there are five $\nu\in Y_\vdom$ that fit in $\lambda$,
giving rise to seven semi-standard Young tableaux with ballot
sequences for reverse row words:
$$
\setlength{\unitlength}{10pt}
\newcommand{\Yd}[1]{\begin{picture}(2,4.5)
\put(0,0){\framebox(1,4){}}\put(0,2){\framebox(2,2){}}
\put(0,1){\line(1,0){1}}\put(0,3){\line(1,0){2}} #1 \end{picture}}
\newcommand{\scribble}[2]{\multiput(#1,#2)(0,.33){6}{\line(3,1){1}}}
\newcommand{\tabl}[3]{\put(#1,#2){\makebox(1,1){\small #3}}}
\Yd{\tabl131\tabl031\tabl122\tabl022\tabl013\tabl004}
\qquad
\Yd{\scribble02\tabl131\tabl122\tabl011\tabl002}
\qquad
\Yd{\scribble02\tabl131\tabl122\tabl011\tabl003}
\qquad
\Yd{\scribble02\tabl131\tabl122\tabl013\tabl004}
\qquad
\Yd{\scribble02\scribble12\tabl011\tabl002}
\qquad
\Yd{\scribble02\scribble00\tabl131\tabl122}
\qquad
\Yd{\scribble02\scribble12\scribble00}
$$
Two tableaux have content $\mu=\omega_2$, and therefore $V(\omega_2)$
occurs with multiplicity two in $W(\omega_2+\omega_4)$.  We leave it
to the enterprising reader to check that in general,
$$
W(a\omega_2+b\omega_4) = \sum_\mu m_\mu V(\mu),
$$
where the sum is over all 
$\mu=c_1\omega_1+c_2\omega_2+c_3\omega_3+c_4\omega_4$ such that
$\mu\subseteq a\omega_2+b\omega_4$ and $c_1=c_3\leq a$, and with
multiplicities $m_\mu$ given by
$$
m_\mu = 1+\min(c_2,\,a-c_3,\,b-c_3-c_4,\,a+b-c_1-c_2-c_3-c_4).
$$
\end{example}


The same techniques used in in the previous example can be used here.
One first identifies the minimal subset of roots $\beta_{k,l}$ such
that $x_{\beta_{k,l}}^-\otimes (t-1).w_\lambda\ne 0$. A simple
counting argument then gives us the maximal possible value for
$m_{\lambda,\mu}$. Then one can prove analogues of Lemma \ref{reduce}
in exactly the same way, to give an upper bound on $m_{\lambda,\mu}$.
When $a=b=1$ this upper bound is precisely the multiplicity $m_\mu$
given above.  In particular, when $\mu=\omega_2$, the bound
$m_{\lambda,\mu}\leq 2$ arises because the difference $\lambda-\mu$
can be written either as $\beta_{1,2}+\beta_{3,4}$ or as
$\beta_{1,3}+\beta_{2,4}$.

}


\begin{thebibliography}{99}

\bibitem{CP}
Chari, V.; Pressley, A.
Minimal affinizations of representations of quantum groups: the
simply-laced case.
{\em J. Algebra} {\bf 184} (1996), no. 1, 1--30.

\bibitem{CPweyl}
Chari, V.; Pressley, A.
Weyl modules for classical and quantum affine algebras.
preprint, math.QA/0004174.

\bibitem{C}
Chari, V.
On the fermionc formula and the Kirillov--Reshetikhin conjecture.
preprint, math.QA/0006090

\bibitem{Fulton}
Fulton, W.
{\em Young Tableaux}.  London Mathematical Society Student Texts \#35,
Cambridge University Press, Cambridge, 1997.

\bibitem{FH}
Fulton, W.; Harris, J.
{\em Representation Theory: A First Course}.
Graduate Texts in Mathematics \#129, Springer--Verlag, New York, 1991.

\bibitem{HKOTY}
Hatayama, G.; Kuniba, A.; Okado, M.; Takagi, T.; Yamada, Y.
Remarks on fermionic formula.
{\em Recent developments in quantum affine algebras and related topics
(Raleigh, NC, 1998)},
Contemp. Math. {\bf 248}, 243--291,  

\bibitem{H}
Humphreys, J.
{\em Introduction to Lie algebras and representation theory}.
Graduate Texts in Mathematics \#9, Springer--Verlag, New York, 1972.

\bibitem{KR}
Kirillov, A. N.; Reshetikhin, N. Yu.
Representations of Yangians and multiplicities of occurrence of the
irreducible components of the tensor product of representations of simple
Lie algebras.
(Russian) {\em Zap. Nauchn. Sem. Leningrad. Otdel. Mat. Inst.
Steklov. (LOMI)} {\bf 160} (1987), transl. in
{\em J. Soviet Math.} {\bf 52} (1990), 3156--3164.

\bibitem{KT}
Koike, K.; Terada, I.
Young-diagrammatic methods for the representation theory of the
classical groups of type $B\sb n,\;C\sb n,\;D\sb n$.
{\em J. Algebra} {\bf 107} (1987), no. 2, 466--511.

\bibitem{kleber}
Kleber, M.
Combinatorial Structure of Finite Dimensional Representations
of Yangians: the Simply-Laced Case.
{\em International Mathematics Research Notices (IMRN)} 1997 \#4, 187--201.

\bibitem{embed}
Kleber, M.
Embeddings of Schur functions into types $B/C/D$.
preprint, math.CO/0009199.

\bibitem{L}
Lusztig, G.
{\em Introduction to quantum groups}.
Progress in Mathematics {\bf 110}, Birkh\"auser, Boston, 1993.

\bibitem{Macdonald}
Macdonald, I. G.
{\em Symmetric Functions and Hall Polynomials,} second edition.
Oxford University Press, Oxford, 1995.

\bibitem{V}
Varagnolo, M.
Quiver Varieties and Yangians.  e-print math.QA/0005277.

\end{thebibliography}
\end{document}